\documentclass[12pt]{article}
\usepackage{amssymb}
\usepackage{amsmath}
\usepackage{graphicx}
\usepackage{subcaption}
\usepackage{float}
\begin{document}

\begin{center}
\textbf{The canonical projection associated with}

\textbf{a} \textbf{mixed possibly infinite iterated function system}

\bigskip

by

\bigskip

\textit{Bogdan Cristian ANGHELINA}, \textit{Radu MICULESCU}, and \textit{%
Alexandru MIHAIL}
\end{center}

\bigskip

\textbf{Abstract}. {\small This paper provides an alternative description
for the fixed points of the fractal operator associated with a mixed
possibly infinite iterated function system via a canonical projection type
function. Some visual aspects of our results are presented.}

\bigskip



\bigskip

\textbf{1. Introduction}

\bigskip

An iterated function system is a pair $((X,d),(f_{i})_{i\in I})=\mathcal{S}$%
, where $(X,d)$ is a complete metric space, $I$ is a finite set and $%
f_{i}:X\rightarrow X$, $i\in I$, are Banach contractions. With such a system
we can associate the fractal operator $F_{\mathcal{S}}:P_{cp}(X)\rightarrow
P_{cp}(X)$, given by $F_{\mathcal{S}}(K)=\underset{i\in I}{\cup }f_{i}(K)$,
for all $K\in P_{cp}(X)=\{C\mid C$ is a non empty compact subset of $X\}$.
Moreover, if $\mathcal{S}$ is endowed with probabilities, we can also
consider the Markov operator acting on a certain set of probability Borel
measures. Via the Banach contraction principle, J. Hutchinson (see [5])
proved that $F_{\mathcal{S}}$ has a unique fixed point $A_{\mathcal{S}}$
(which is called the attractor of $\mathcal{S}$) and that the Markov
operator also has a unique fixed point which is called the associated
Hutchinson measure. The sets that can be represented as attractors of
iterated function systems are called Hutchinson-Barnsley fractals.

The canonical projection associated with $\mathcal{S}$ is a continuous
surjection $\pi $ from the code space $\Lambda (I)$ onto the attractor $A_{%
\mathcal{S}}$. It is also called \textquotedblleft the coding
map\textquotedblright\ (see [12] and [13]), \textquotedblleft the address
map\textquotedblright\ (see [9]) or \textquotedblleft the coordinate
map\textquotedblright\ (see [8]).

The importance of this projection is emphasized by the following facts:

a) It provides, via the formula $\pi (\Lambda (I))=A_{\mathcal{S}}$, an
alternative description of the attractor of $\mathcal{S}$ and, consequently,
it is an important device in the topological study of Hutchinson-Barnsley
fractals.

b) It is involved in the alternative presentation of the Hutchinson measure
associated with $\mathcal{S}$ as the push-forward measure of the Bernoulli
measure on $\Lambda (I)$ through $\pi $.

c) It was a source of inspiration for the concept of topological
self-similar set introduced by A. Kameyama (see [7]).

The above mentioned theory (which was initiated by J. Hutchinson and
developed by M. Barnsley) has been extended in different directions. Two of
them are of special interest from the point of view of the present paper:

i) the direction concentrating on systems involving not necessarily finite
families of functions (see [2], [3], [4], [10], [12], [14] and [15]);

ii) the direction focusing on systems involving functions from larger
classes of contractions (see [11] which, together with [6], is an excellent
survey on iterated function systems, and the references therein).

We emphasize that the corresponding fractal operator is a Picard operator
for all the above mentioned generalizations of the concept of iterated
function system.

In the recent work [1], we integrated the previous directions by considering
possibly infinite iterated function systems enriched with orbital possibly
infinite iterated function systems (called mixed possibly infinite iterated
function systems) and we proved that the fractal operator associated with
such a system is weakly Picard. Its fixed points are called attractors.

The present paper develops a canonical projection type theory for mixed
possibly infinite iterated function systems in order to obtain an
alternative description for the attractors. Finally some visual aspects
concerning our results are presented.

\bigskip

\textbf{2. Preliminaries}

\bigskip

\textbf{The generalized pigeonhole principle}

\bigskip

In the sequel we will use the following form of the generalized pigeonhole
principle: If $N$ objects are placed in $k$ boxes, then at least one box
contains at least $[\frac{N-1}{k}]+1$.

\bigskip

\textbf{Basic notation}

\bigskip

By $\mathbb{N}$ we mean the set $\{1,2,...\}$.

For a function $f:X\rightarrow X$ and $n\in \mathbb{N}$, by $f^{[n]}$ we
mean the composition of $f$ by itself $n$ times and by $f^{[0]}$ we mean $%
Id_{X}$.

For a Lipschitz function $f:X\rightarrow Y$, where $(X,d)$ and $(Y,\rho )$
are metric spaces, by $lip(f)\,$we denote the Lipschitz constant of $f$.

For a metric space $(X,d)$ and $A,B\subseteq X$ we shall use the following
notation:

- $\{A\subseteq X\mid A\neq \emptyset $ and $A$ is bounded$\}\overset{not}{=}%
P_{b}(X)$

- $\{A\subseteq X\mid A\neq \emptyset $ and $A$ is closed and bounded$\}%
\overset{not}{=}P_{b,cl}(X)$

- $\{A\subseteq X\mid A\neq \emptyset $ and $A$ is compact$\}\overset{not}{=}%
P_{cp}(X)$

- $\underset{a\in A}{\inf }$ $d(x,a)\overset{not}{=}d(x,A)$

- $\underset{a\in A}{\sup }$ $d(a,B)\overset{not}{=}D(A,B)$

- $\underset{x,y\in A}{\sup }d(x,y)\overset{not}{=}diam(A)$.

\bigskip

\textbf{The code space}

\bigskip

Given a set $I$ and $n\in \mathbb{N}$, we consider:

- 
\begin{equation*}
I^{\mathbb{N}}\overset{not}{=}\Lambda (I)
\end{equation*}

-%
\begin{equation*}
I^{\{1,...,n\}}\overset{not}{=}\Lambda _{n}(I)\text{.}
\end{equation*}

So:

- the elements of $\Lambda (I)$, which is called a code space, are written
as infinite words $\alpha =\alpha _{1}\alpha _{2}...\alpha _{n}\alpha
_{n+1}...$ with letters from $I$

- the elements of $\Lambda _{n}(I)$ are written as words $\alpha =\alpha
_{1}\alpha _{2}...\alpha _{n}$ having $n$ letters from $I$ and $n$, which is
denoted by $\left\vert \alpha \right\vert $, is called the length of $\alpha 
$.

In the sequel we will use the following notation:%
\begin{equation*}
\underset{n\in \mathbb{N\cup }\{0\}}{\cup }\Lambda _{n}(I)\overset{not}{=}%
\Lambda ^{\ast }(I)\text{,}
\end{equation*}%
where $\Lambda _{0}(I)=\{\lambda \}$ and $\lambda $ designates the empty
word.

If $\alpha =\alpha _{1}\alpha _{2}...\alpha _{n}\alpha _{n+1}....\in \Lambda
(I)$ or if $\alpha =\alpha _{1}\alpha _{2}...\alpha _{n}\in \Lambda _{n}(B)$
and $m,n\in \mathbb{N}$, $n\geq m$, then we will use the following notation: 
\begin{equation*}
\alpha _{1}\alpha _{2}...\alpha _{m}\overset{not}{=}[\alpha ]_{m}\text{.}
\end{equation*}

By the concatenation of the words $\alpha $ and $\beta $, where $\alpha
=\alpha _{1}\alpha _{2}...\alpha _{n}\in \Lambda _{n}(I)\,$and $\beta =\beta
_{1}\beta _{2}...\beta _{m}\beta _{m+1}...\in \Lambda (I)$, we mean$\ $the
infinite word $\alpha _{1}...\alpha _{n}\beta _{1}...\beta _{m}\beta
_{m+1}...$ which is denoted by $\alpha \beta $.

For $i\in I$, we introduce the function $\tau _{i}:\Lambda (I)\rightarrow
\Lambda (I)$,\ given by 
\begin{equation*}
\tau _{i}(\alpha )=i\alpha \text{,}
\end{equation*}
for all $\alpha \in \Lambda (I)$.

$\Lambda (I)$ can be endowed with the metric $d_{\Lambda (I)}$, called the
Baire metric, described by: 
\begin{equation*}
d_{\Lambda }(\alpha ,\beta )=\{%
\begin{array}{cc}
0\text{,} & \text{if }\alpha =\beta \\ 
\frac{1}{2^{\min \{k\in \mathbb{N}\mid \alpha _{k}\neq \beta _{k}\}}}\text{,}
& \text{if }\alpha \neq \beta%
\end{array}%
\text{,}
\end{equation*}%
for all $\alpha =\alpha _{1}\alpha _{2}...\alpha _{n}\alpha _{n+1}...\in
\Lambda (I)$ and $\beta =\beta _{1}\beta _{2}...\beta _{n}\beta _{n+1}...\in
\Lambda (I)$.

Given $f_{i}:X\rightarrow X$, $i\in I$, and $\alpha =\alpha _{1}\alpha
_{2}...\alpha _{n}\in \Lambda _{n}(I)$, the following notation will be used
in the sequel:%
\begin{equation*}
f_{\alpha _{1}}\circ f_{\alpha _{2}}\circ ....\circ f_{\alpha _{n}}\overset{%
not}{=}f_{\alpha }
\end{equation*}%
and%
\begin{equation*}
Id_{X}\overset{not}{=}f_{\lambda }\text{.}
\end{equation*}

In particular if $I=\{i\}$, we have $f_{\underset{n\text{ times}}{i...i}%
}=f_{i}^{[n]}$ for every $n\in \mathbb{N}$.

\bigskip

\textbf{The Hausdorff-Pompeiu metric}

\bigskip

\textbf{Remark 2.1.}

\textit{i)} \textit{Given a metric space }$(X,d)$ \textit{and }$A$\textit{\
and }$B$\textit{\ subsets of }$X$\textit{, we have}%
\begin{equation*}
D(A,B)=D(\overline{A},\overline{B})\text{.}
\end{equation*}

\textit{ii)} \textit{Given a metric space }$(X,d)$ \textit{and a set }$I$%
\textit{, we have }%
\begin{equation*}
D(\underset{i\in I}{\cup }A_{i},\underset{i\in I}{\cup }B_{i})\leq \underset{%
i\in I}{\sup }\text{ }D(A_{i},B_{i})\text{,}
\end{equation*}%
\textit{for every }$A_{i}$ \textit{and }$B_{i}$ \textit{subsets of }$X$%
\textit{.}

\bigskip

\textbf{Proposition 2.2.} \textit{Given a metric space }$(X,d)$\textit{, the
function} $h^{\ast }:P_{b}(X)\times P_{b}(X)\rightarrow \lbrack 0,\infty )$%
\textit{, described by\ }%
\begin{equation*}
h^{\ast }(A,B)=\max \{D(A,B),D(B,A)\}\text{,}
\end{equation*}%
\textit{for every }$A,B\in P_{b}(X)$\textit{, has the following properties:}

\textit{i) }%
\begin{equation*}
h^{\ast }(\{x\},\{y\})=d(x,y)\text{,}
\end{equation*}%
\textit{for every }$x,y\in X$\textit{;}

\textit{ii) }%
\begin{equation*}
h^{\ast }(A,B)=h^{\ast }(\overline{A},\overline{B})\text{,}
\end{equation*}%
\textit{\ for every }$A,B\in P_{b}(X)$\textit{;}

\textit{iii) }%
\begin{equation*}
h^{\ast }(\underset{i\in I}{\cup }A_{i},\underset{i\in I}{\cup }B_{i})\leq 
\underset{i\in I}{\sup }\text{ }h^{\ast }(A_{i},B_{i})\text{,}
\end{equation*}%
\textit{for all families }$\{A_{i}\}_{i\in I}$\textit{\ and }$%
\{B_{i}\}_{i\in I}$\textit{\ of elements of }$P_{b}(X)$\textit{\ such that }$%
\underset{i\in I}{\cup }A_{i}\in P_{b}(X)$\textit{\ and }$\underset{i\in I}{%
\cup }B_{i}\in P_{b}(X)$\textit{.}

\bigskip 

Let us recall two well-known results.

\bigskip

\textbf{Definition 2.3.} \textit{Given a metric space }$(X,d)$\textit{, the
restriction of }$h^{\ast }$\textit{\ to }$P_{b,cl}(X)\times P_{b,cl}(X)$%
\textit{\ is a metric, denoted by }$h$\textit{, called the Hausdorff-Pompeiu
metric.}

\bigskip

\textbf{Proposition 2.4.} \textit{Given a complete metric space }$(X,d)$%
\textit{, the metric space }$(P_{b,cl}(X),h)$ \textit{is complete and if }$%
(A_{n})_{n\in \mathbb{N}}\subseteq P_{b,cl}(X)$\textit{\ is Cauchy, then}%
\begin{equation*}
\underset{n\rightarrow \infty }{\lim }A_{n}=\{x\in X\mid \text{\textit{there
exist an increasing sequence} }(n_{k})_{k\in \mathbb{N}}\subseteq \mathbb{N}
\end{equation*}%
\begin{equation*}
\text{\textit{and} }x_{n_{k}}\in A_{n_{k}}\text{ \textit{for every} }k\in 
\mathbb{N}\text{\textit{\ such that} }\underset{k\rightarrow \infty }{\lim }%
x_{n_{k}}=x\}\text{.}
\end{equation*}

\bigskip

\textbf{Proposition 2.5.} \textit{Given a metric space }$(X,d)$\textit{\ and 
}$(A_{n})_{n\in \mathbb{N}}\subseteq P_{b,cl}(X)$ \textit{such that }$%
A_{n+1}\subseteq A_{n}$ \textit{for every} $n\in \mathbb{N}$\textit{, we have%
}%
\begin{equation*}
\underset{n\rightarrow \infty }{\lim }h(A_{n},\underset{n\in \mathbb{N}}{%
\cap }A_{n})=0\text{.}
\end{equation*}

\bigskip

\textbf{Possibly infinite iterated function systems}

\bigskip

\textbf{Definition 2.6.} \textit{A possibly infinite iterated function
system (for short IIFS) is a pair} $((X,d),(f_{i})_{i\in I})\overset{not}{=}%
\mathcal{S}$, \textit{where }$(X,d)$\textit{\ is a complete metric space and 
}$f_{i}:X\rightarrow X$\textit{, }$i\in I$\textit{, are such that:}

\textit{a)} $f_{i}$ \textit{is continuous for every }$i\in I$;

\textit{b) the family }$(f_{i})_{i\in I}$\textit{\ is bounded, i.e. }%
\begin{equation*}
\underset{i\in I}{\cup }f_{i}(B)\in P_{b}(X)\text{,}
\end{equation*}%
\textit{\ for every }$B\in P_{b}(X)$.

\textit{The function }$F_{\mathcal{S}}:P_{b,cl}(X)\rightarrow P_{b,cl}(X)$%
\textit{, given by}%
\begin{equation*}
F_{\mathcal{S}}(B)=\overline{\underset{i\in I}{\cup }f_{i}(B)}\text{,}
\end{equation*}%
\textit{for every }$B\in P_{b,cl}(X)$\textit{, is called the} \textit{%
fractal operator associated with} $\mathcal{S}$.

\bigskip

\textbf{Remark 2.7.} \textit{In the framework of the above definition, if
the functions }$f_{i}$\textit{\ are Lipschitz, then}%
\begin{equation*}
lip(F_{\mathcal{S}})\leq \underset{i\in I}{\sup }\text{ }lip(f_{i})\text{.}
\end{equation*}

\bigskip

\textbf{Proposition 2.8 }(see Lemma 2.6 from [17])\textbf{.} \textit{For each%
} \textit{IIFS} $\mathcal{S}=((X,d),(f_{i})_{i\in I})$\textit{, we have} 
\begin{equation*}
F_{\mathcal{S}}^{[n]}(B)=\overline{\underset{\alpha \in \Lambda _{n}(I)}{%
\cup }f_{\alpha }(B)}\text{,}
\end{equation*}%
\textit{for every\ }$n\in \mathbb{N}$\textit{\ and every\ }$B\in P_{b,cl}(X)$%
\textit{.}

\bigskip

\textbf{Theorem 2.9 }(see [15])\textbf{.} \textit{For each IIFS} $\mathcal{S}%
=(X,(f_{i})_{i\in I})$\textit{\ such that the functions }$f_{i}$ \textit{are
Lipschitz and} $\underset{i\in I}{\sup }$ $lip(f_{i})<1$\textit{, }$F_{%
\mathcal{S}}$ \textit{is a Banach contraction with respect to }$h$ \textit{%
and its unique fixed point} $A_{\mathcal{S}}$\textit{\ is called the
attractor of} $\mathcal{S}$.

\textit{Moreover, we have:}

\textit{a)} \textit{For every }$\alpha \in \Lambda (I)$\textit{, the set }$%
\underset{n\in \mathbb{N}}{\cap }\overline{f_{[\alpha ]_{n}}(A_{\mathcal{S}})%
}$\textit{\ has just one element, which is denoted by }$a_{\alpha }$\textit{.%
}

\textit{b)} 
\begin{equation*}
\underset{n\rightarrow \infty }{\lim }f_{[\alpha ]_{n}}(x)=a_{\alpha }\text{,%
}
\end{equation*}%
\textit{for every }$x\in X$\textit{\ and every }$\alpha \in \Lambda (I)$%
\textit{.}

\textit{c)} \textit{The continuous (with respect to the Baire metric)
function }$\pi :\Lambda (I)\rightarrow A_{\mathcal{S}}$\textit{, defined by }%
\begin{equation*}
\pi (\alpha )=a_{\alpha }\text{,}
\end{equation*}%
\textit{for every\ }$\alpha \in \Lambda (I)$\textit{, is called the
canonical projection associated with }$\mathcal{S}$ \textit{and it\ has the
following properties:}

\qquad \textit{i)} 
\begin{equation*}
\pi \circ \tau _{i}=f_{i}\circ \pi \text{,}
\end{equation*}%
\textit{\ for every }$i\in I$;

\qquad \textit{ii)\ }%
\begin{equation*}
\overline{\pi (\Lambda (I))}=A_{\mathcal{S}}\text{.}
\end{equation*}

\bigskip

Given a metric space $(X,d)$, $x\in X$, $B\subseteq X$ and a family of
functions $\mathcal{F}=(f_{i})_{i\in I}$\ , where $f_{i}:X\rightarrow X$, we
shall use the notation:

\begin{equation*}
\underset{n\in \mathbb{N\cup }\{0\}}{\cup }\overline{\underset{\alpha \in
\Lambda _{n}(I)}{\cup }f_{\alpha }(B)}\overset{not}{=}\mathcal{O}_{\mathcal{F%
}}(B)
\end{equation*}%
and%
\begin{equation*}
\mathcal{O}_{\mathcal{F}}(\{x\})\overset{not}{=}\mathcal{O}_{\mathcal{F}}(x)%
\text{.}
\end{equation*}

In particular, given an IIFS $\mathcal{S}=((X,d),(f_{i})_{i\in I})$ and $%
B\in P_{b,cl}(X)$, we shall use the notation%
\begin{equation*}
\mathcal{O}_{(f_{i})_{i\in I}}(B)\overset{not}{=}\mathcal{O}_{\mathcal{S}}(B)%
\text{.}
\end{equation*}

\bigskip

\textbf{Definition 2.10.} \textit{An orbital possibly infinite iterated
function system (for short oIIFS) is an IIFS }$\mathcal{S=}%
((X,d),(f_{i})_{i\in I})$ \textit{such that:}

\textit{a) the family }$(f_{i})_{i\in I}$\textit{\ is equi-uniformly
continuous on bounded sets, i.e. for every }$B\in P_{b}(X)$\textit{\ and
every }$\varepsilon >0$\textit{\ there exists }$\delta _{\varepsilon ,B}>0$%
\textit{\ such that for all }$x,y\in B$ \textit{and} $i\in I$\textit{\ the
following implication is valid: }$d(x,y)<\delta _{\varepsilon ,B}\Rightarrow
d(f_{i}(x),f_{i}(y))<\varepsilon $\textit{.}

\textit{b) there exists }$a\in \lbrack 0,1)$ \textit{such that} 
\begin{equation*}
d(f_{i}(y),f_{i}(z))\leq ad(y,z)\text{,}
\end{equation*}%
\textit{for every }$i\in I$\textit{, }$x\in X$\textit{\ and} $y,z\in 
\mathcal{O}_{\mathcal{(}f_{i}\mathcal{)}_{i\in I}}(x)$.

\bigskip

\textbf{Definition 2.11.} \textit{A mixed possibly infinite iterated
function system (briefly mIIFS) is an IIFS }$\mathcal{S=}((X,d),(f_{i})_{i%
\in I\cup J})$, \textit{where }$I$\textit{\ and }$J$\textit{\ are disjoint
sets, such that:}

\textit{a)}%
\begin{equation*}
lip(f_{i})\leq 1\text{,}
\end{equation*}%
\textit{for every }$i\in I\cup J$;

\textit{b) there exists }$a\in \lbrack 0,1)$\textit{\ having the following
properties:}

\qquad \textit{b 1)\ }%
\begin{equation*}
lip(f_{i})\leq a\text{,}
\end{equation*}%
\textit{for every }$i\in I$;

\qquad \textit{b 2) }%
\begin{equation*}
lip(f_{i\mid \overline{\mathcal{O}_{\mathcal{(}f_{i}\mathcal{)}_{i\in J}}(x)}%
})\leq a\text{,}
\end{equation*}%
\textit{for every }$x\in X$\textit{\ and every }$i\in J$.

\bigskip

\textbf{Proposition 2.12 }(see [1])\textbf{.} \textit{For each mIFS }$%
\mathcal{S=}((X,d),(f_{i})_{i\in I\cup J})$ \textit{and every} $B\in
P_{b,cl}(X)$ \textit{there exists} $A_{B}\in P_{b,cl}(X)$ \textit{such that }%
\begin{equation*}
F_{\mathcal{S}}(A_{B})=A_{B}
\end{equation*}%
\textit{and} 
\begin{equation*}
\underset{n\rightarrow \infty }{\lim }h(F_{\mathcal{S}}^{[n]}(B),A_{B})=0,
\end{equation*}%
\textit{\ i.e.} $F_{\mathcal{S}}$ \textit{is weakly Picard.}

\bigskip

Note that, in the framework of the above definition, $((X,d),(f_{i})_{i\in
I})\overset{not}{=}\mathcal{S}_{I}$ and $((X,d),(f_{\omega i\theta })_{i\in
I,\omega ,\theta \in \Lambda ^{\ast }(J)})\overset{not}{=}\mathfrak{S}$ are
IIFSs (see Remark 2.9 and Lemma 3.5 from [1]) and $((X,d),(f_{i})_{i\in J})%
\overset{not}{=}\mathcal{S}_{J}$ is an oIIFS.

\bigskip

Given an mIIFS $\mathcal{S=}((X,d),(f_{i})_{i\in I\cup J})$, $x\in X$, $B\in
P_{b,cl}(X)$ and $\alpha =\alpha _{1}\alpha _{2}...\alpha _{n}...\in \Lambda
(I\cup J)$ or $\alpha =\alpha _{1}\alpha _{2}...\alpha _{n}\in \Lambda
^{\ast }(I\cup J)$, we shall use the following notation:

-%
\begin{equation*}
A_{\{x\}}\overset{not}{=}A_{x}
\end{equation*}

- 
\begin{equation*}
\mathcal{O}_{\mathcal{S}_{J}}(B)\overset{not}{=}\mathcal{O}_{J}(B)\overset{%
\text{Lemma 3.36 from [16]}}{\in }P_{b}(X)
\end{equation*}

-%
\begin{equation*}
\underset{x\in B}{\sup }\max \{diam(\{x\}\cup F_{\mathcal{S}}(\{x\}),diam(%
\mathcal{O}_{\mathfrak{S}}(x)),diam(\mathcal{O}_{J}(x))\}\overset{not}{=}%
N_{B}\in \mathbb{R}
\end{equation*}

-%
\begin{equation*}
card(\{l\in \mathbb{N}\mid \alpha _{l}\in I\})\overset{not}{=}n_{I}(\alpha )%
\text{.}
\end{equation*}

- $\Lambda _{1}(I\cup J)$ denotes the set of finite words with letters from $%
I\cup J$ ending with a letter from $I$

- $\Lambda _{2}(I\cup J)$ denotes the set of finite words with letters from $%
I\cup J$ starting with a letter from $I$

- $\Lambda _{3}(I\cup J)=\Lambda _{1}(I\cup J)\cap \Lambda _{2}(I\cup J)$
denotes the set of finite words with letters from $I\cup J$ starting and
ending with letters from $I$

- $\Sigma _{0}(I,J)$ denotes the set of finite words with letters from $%
I\cup J\cup \Lambda (J)$ having the form 
\begin{equation*}
\beta _{0}\gamma _{1}\beta _{1}...\gamma _{n}\beta _{n}\text{,}
\end{equation*}%
with $n\in \mathbb{N}$, where%
\begin{equation*}
\beta _{0}\in \{\lambda \}\cup \Lambda _{1}(I\cup J)\text{,}
\end{equation*}%
\begin{equation*}
\beta _{n}\in \{\lambda \}\cup \Lambda _{2}(I\cup J)\text{,}
\end{equation*}%
\begin{equation*}
\gamma _{k}\in \Lambda (J)
\end{equation*}%
for every $k\in \{1,...,n\}$ and%
\begin{equation*}
\beta _{k}\in \Lambda _{3}(I\cup J)\text{,}
\end{equation*}%
for every $k\in \{1,...,n-1\}$ if $n\geq 2$

- 
\begin{equation*}
\Sigma _{1}(I,J)=\{\alpha \in \Lambda (I\cup J)\mid n_{I}(\alpha )=\infty \}
\end{equation*}

- 
\begin{equation*}
\Sigma (I,J)=\Sigma _{0}(I,J)\cup \Sigma _{1}(I,J)\text{.}
\end{equation*}

Additionally, we consider the function $g:X\rightarrow P_{b,cl}(X)$, given by%
\begin{equation*}
g(x)=A_{x}\text{,}
\end{equation*}%
for every $x\in X$.

\bigskip

\textbf{Lemma 2.13.}\textit{\ Given an mIIFS} $\mathcal{S=}%
((X,d),(f_{i})_{i\in I\cup J})$\textit{, we have}%
\begin{equation*}
lip(g)\leq 1\text{.}
\end{equation*}

\textit{Proof}. Indeed, we have%
\begin{equation*}
h(g(x),g(y))=h(A_{x},A_{y})\leq 
\end{equation*}%
\begin{equation*}
\leq h(A_{x},F_{\mathcal{S}}^{[n]}(\{x\}))+h(F_{\mathcal{S}}^{[n]}(\{x\}),F_{%
\mathcal{S}}^{[n]}(\{y\}))+h(F_{\mathcal{S}}^{[n]}(\{y\}),A_{y})\leq 
\end{equation*}%
\begin{equation*}
\overset{\text{Remark 2.7 \& Proposition 2.2, i)}}{\leq }h(A_{x},F_{\mathcal{%
S}}^{[n]}(\{x\}))+h(F_{\mathcal{S}}^{[n]}(\{y\}),A_{y})+d(x,y)\text{,}
\end{equation*}%
for all $n\in \mathbb{N}$, so, by passing to the limit as $n\rightarrow
\infty $, we get%
\begin{equation*}
h(g(x),g(y))\leq d(x,y)\text{,}
\end{equation*}%
for all $x,y\in X$. $\square $

\bigskip

\textbf{Lemma 2.14.} \textit{Given an mIIFS} $\mathcal{S=}%
((X,d),(f_{i})_{i\in I\cup J})$\textit{, we have}%
\begin{equation*}
A_{f_{i}(x)}\subseteq A_{x}\text{,}
\end{equation*}%
\textit{for all }$x\in X$\textit{\ and }$i\in I\cup J$\textit{.}

\textit{Proof}. We have%
\begin{equation*}
A_{f_{i}(x)}\overset{\text{Proposition 2.12}}{=}\underset{n\rightarrow
\infty }{\lim }F_{\mathcal{S}}^{[n]}(\{f_{i}(x)\})\overset{\text{Proposition
2.4}}{=}\{y\in X\mid \text{there exist }
\end{equation*}%
\begin{equation*}
\text{a strictly increasing sequence }(n_{k})_{k\in \mathbb{N}}\text{ of
natural numbers and}
\end{equation*}%
\begin{equation*}
y_{n_{k}}\in F_{\mathcal{S}}^{[n_{k}]}(\{f_{i}(x)\})\text{ for every }k\in 
\mathbb{N}\text{ such that }y=\underset{k\rightarrow \infty }{\lim }%
y_{n_{k}}\}\subseteq
\end{equation*}%
\begin{equation*}
\overset{F_{\mathcal{S}}^{[n_{k}]}(\{f_{i}(x)\})\subseteq F_{\mathcal{S}%
}^{[n_{k}+1]}(\{x\})}{\subseteq }\{y\in X\mid \text{there exist a strictly
increasing sequence}
\end{equation*}%
\begin{equation*}
(m_{k})_{k\in \mathbb{N}}\text{ of natural numbers and }y_{m_{k}}\in F_{%
\mathcal{S}}^{[m_{k}]}(\{x\})\text{ for every }k\in \mathbb{N}
\end{equation*}%
\begin{equation*}
\text{such that }y=\underset{k\rightarrow \infty }{\lim }y_{m_{k}}\}=A_{x}%
\text{,}
\end{equation*}%
for all $x\in X$ and $i\in I\cup J$. $\square $

\bigskip

\textbf{3. The main results}

\bigskip

Our first result is a counterpart of Theorem 2.9, b).

\bigskip

\textbf{Proposition 3.1.} \textit{For each mIIFS} $\mathcal{S}%
=((X,d),(f_{i})_{i\in I\cup J})$, $\alpha \in \Lambda (I\cup J)$ \textit{and 
}$x\in X$\textit{,\ the sequence} $(f_{[\alpha ]_{n}}(x))_{n\in \mathbb{N}}$ 
\textit{is convergent.}

\textit{Proof}. We divide the proof into two cases:

a) $\alpha \in \Sigma _{1}(I,J)$, i.e. $n_{I}(\alpha )$ is infinite;

b) $\alpha \in \Lambda (I\cup J)\setminus \Sigma _{1}(I,J)$, i.e. $%
n_{I}(\alpha )$ is finite.

In the first case, we can find a strictly increasing sequence $(n_{k})_{k\in 
\mathbb{N}}$ of natural numbers such that%
\begin{equation*}
\alpha _{n_{k}}\in I\text{,}
\end{equation*}%
for every $k\in \mathbb{N}$ and 
\begin{equation*}
\alpha _{n}\notin I\text{, i.e. }\alpha _{n}\in J\text{,}
\end{equation*}%
for every $n\in \mathbb{N\smallsetminus \{}n_{k}\mid k\in \mathbb{N\}}$.

\textbf{Claim}. The sequence $(f_{[\alpha ]_{n}}(x))_{n\in \mathbb{N}}$ is
Cauchy.

\textit{Justification of the claim}. Adopting the notation%
\begin{equation*}
\alpha _{n_{k-1}+1}\alpha _{n_{k-1}+2}...\alpha _{n_{k}}\overset{not}{=}%
\beta _{k}\text{,}
\end{equation*}%
we have%
\begin{equation*}
d(f_{[\alpha ]_{n_{k}}}(x),f_{[\alpha ]_{n_{k+1}}}(x))=d(f_{[\alpha
]_{n_{k}}}(x),f_{[\alpha ]_{n_{k}}\beta _{k+1}}(x))\overset{\text{Definition
2.11, b 1) }}{\leq }
\end{equation*}%
\begin{equation*}
\leq a^{k}d(x,f_{\beta _{k+1}}(x))\leq a^{k}diam(O_{\mathfrak{S}}(x))\text{,}
\end{equation*}%
so%
\begin{equation*}
d(f_{[\alpha ]_{n_{k}}}(x),f_{[\alpha ]_{n_{k+p}}}(x))\leq
(a^{k}+a^{k+1}+...+a^{k+p-1})diam(O_{\mathfrak{S}}(x))\leq
\end{equation*}%
\begin{equation}
\leq \frac{a^{k}}{1-a}diam(O_{\mathfrak{S}}(x))\text{,}  \tag{1}
\end{equation}%
for every $k,p\in \mathbb{N}$.

For $k,m,n\in \mathbb{N}$ such that $n_{k}\leq n\leq m$, there exist $s,t\in 
\mathbb{N}$, $k\leq s\leq t$ having the property that $n_{s}\leq n<n_{s+1}$
and $n_{t}\leq m<n_{t+1}$. Hence we get%
\begin{equation}
d(f_{[\alpha ]_{n}}(x),f_{[\alpha ]_{n_{s}}}(x))\leq a^{k}diam(O_{J}(x)) 
\tag{2}
\end{equation}%
and%
\begin{equation}
d(f_{[\alpha ]_{m}}(x),f_{[\alpha ]_{n_{t}}}(x))\leq a^{k}diam(O_{J}(x))%
\text{.}  \tag{3}
\end{equation}%
Consequently%
\begin{equation*}
d(f_{[\alpha ]_{m}}(x),f_{[\alpha ]_{n}}(x))\leq d(f_{[\alpha
]_{m}}(x),f_{[\alpha ]_{n_{t}}}(x))+d(f_{[\alpha ]_{n_{s}}}(x),f_{[\alpha
]_{n_{t}}}(x))+
\end{equation*}%
\begin{equation}
+d(f_{[\alpha ]_{n_{s}}}(x),f_{[\alpha ]_{n}}(x))\overset{\text{(1), (2) \&
(3)}}{\leq }2a^{k}diam(O_{J}(x))+\frac{a^{k}}{1-a}diam(O_{\mathfrak{S}}(x))%
\text{.}  \tag{4}
\end{equation}

The last inequality yields claim's validity.

In view of the Claim we infer that, in the first case, $(f_{[\alpha
]_{n}}(x))_{k\in \mathbb{N}}$ is convergent.

In the second case, there exist $\beta \in \Lambda ^{\ast }(I\cup J)$ and $%
\gamma \in \Lambda (J)$ such that $\alpha =\beta \gamma $ and we get%
\begin{equation*}
(f_{[\alpha ]_{\left\vert \beta \right\vert +n}}(x),f_{[\alpha ]_{\left\vert
\beta \right\vert +n+p}}(x))\overset{\text{Definition 2.11, a) }}{\leq }%
d(f_{[\gamma ]_{n}}(x),f_{[\gamma ]_{n+p}}(x))\leq
\end{equation*}%
\begin{equation}
\overset{\text{Definition 2.11, b 2) }}{\leq }a^{n}diam(O_{J}(x))\text{,} 
\tag{5}
\end{equation}%
for every $n,p\in \mathbb{N}$. The previous inequality ensures that $%
(f_{[\alpha ]_{n}}(x))_{n\in \mathbb{N}}$ is Cauchy, so it is also
convergent.

Now the proof is complete. $\square $

\bigskip

For an mIIFS $\mathcal{S}=((X,d),(f_{i})_{i\in I\cup J})$ and $\alpha \in
\Lambda (I\cup J)$, based on Proposition 3.1, we can consider the function $%
a_{\alpha }:X\rightarrow X$, given by

\begin{equation*}
a_{\alpha }(x)=\underset{n\rightarrow \infty }{\lim }f_{[\alpha ]_{n}}(x)%
\text{,}
\end{equation*}%
for every $x\in X$.

\bigskip

A closer look at the inequalities $(4)$ and $(5)$ from the proof of the
above Proposition leads to the conclusion that the convergence of the
sequence $(f_{[\alpha ]_{n}}(x))_{n\in \mathbb{N}}$ is uniform with respect
to $x$ in a bounded subset of $X$. More precisely, we can state the
following:

\bigskip

\textbf{Corollary 3.2.} \textit{For each mIIFS} $\mathcal{S}%
=((X,d),(f_{i})_{i\in I\cup J})$, $\alpha \in \Lambda (I\cup J)$ \textit{and 
}$B\in P_{b}(X)$\textit{, we have}%
\begin{equation*}
\underset{n\rightarrow \infty }{\lim }\text{ }\underset{x\in B}{\sup }\text{ 
}d(f_{[\alpha ]_{n}}(x),a_{\alpha }(x))=0\text{.}
\end{equation*}

\bigskip

\textbf{Proposition 3.3.} \textit{For an mIIFS} $\mathcal{S}%
=((X,d),(f_{i})_{i\in I\cup J})$, $\alpha \in \Lambda (I\cup J)$ \textit{and 
}$B\in P_{b,cl}(X)$\textit{, we have}%
\begin{equation*}
\overline{a_{\alpha }(B)}\subseteq A_{B}\text{.}
\end{equation*}

\textit{Proof}. For each $B\in P_{b,cl}(X)$, $\alpha \in \Lambda (I\cup J)$, 
$x\in B$ and $n\in \mathbb{N}$, we have%
\begin{equation*}
f_{[\alpha ]_{n}}(x)\overset{\text{Proposition 2.8}}{\in }F_{\mathcal{S}%
}^{[n]}(B)\text{,}
\end{equation*}%
so, since $\underset{n\rightarrow \infty }{\lim }F_{\mathcal{S}}^{[n]}(B)%
\overset{\text{Proposition 2.12}}{=}A_{B}$, via Proposition 2.4, we infer
that%
\begin{equation*}
a_{\alpha }(x)\in A_{B}\text{.}
\end{equation*}

Therefore%
\begin{equation*}
a_{\alpha }(B)\subseteq A_{B}\text{,}
\end{equation*}%
so%
\begin{equation*}
\overline{a_{\alpha }(B)}\subseteq A_{B}\text{. }\square
\end{equation*}

\bigskip

\textbf{Proposition 3.4.} \textit{For each mIIFS} $\mathcal{S}%
=((X,d),(f_{i})_{i\in I\cup J})$ \textit{and} $\alpha \in \Lambda (I\cup J)$%
\textit{, we have:}

\textit{a)}%
\begin{equation*}
lip(a_{\alpha })\leq 1\text{;}
\end{equation*}

\textit{b) }$a_{\alpha }$\textit{\ is constant for every }$\alpha \in \Sigma
_{1}(I,J)$\textit{.}

\textit{Proof}. a) Note that%
\begin{equation}
d(f_{[\alpha ]_{n}}(x),f_{[\alpha ]_{n}}(y))\overset{\text{Definition 2.11,
a) \& b 1)}}{\leq }a^{n_{I}([\alpha ]_{n})}d(x,y)\leq d(x,y)\text{,}  \tag{1}
\end{equation}%
for all $n\in \mathbb{N}$ and $x,y\in X$.

Via $(1)$, by passing to the limit as $n$ goes to $\infty $, we infer that%
\begin{equation*}
d(a_{\alpha }(x),a_{\alpha }(y))\leq d(x,y)\text{,}
\end{equation*}%
for all $x,y\in X$, so $lip(a_{\alpha })\leq 1$.

b) If $n_{I}(\alpha )=\infty $, then $\underset{n\rightarrow \infty }{\lim }%
n_{I}([\alpha ]_{n})=\infty $, and, by passing again to the limit as $n$
goes to $\infty $ in $(1)$, we conclude that%
\begin{equation*}
a_{\alpha }(x)=a_{\alpha }(y)\text{,}
\end{equation*}%
for all $x,y\in X$, so $a_{\alpha }$\textit{\ }is constant. $\square $

\bigskip

\textbf{Proposition 3.5.} \textit{Given an mIIFS} $\mathcal{S}%
=((X,d),(f_{i})_{i\in I\cup J})$, \textit{we have}%
\begin{equation*}
a_{i\alpha }=f_{i}\circ a_{\alpha }\text{,}
\end{equation*}%
\textit{for all} $i\in I\cup J$\textit{\ and }$\alpha \in \Lambda (I\cup J)$%
\textit{.}

\textit{Proof}. Indeed, we have%
\begin{equation*}
f_{i}(a_{\alpha }(x))=f_{i}(\underset{n\rightarrow \infty }{\lim }f_{[\alpha
]_{n}}(x))\overset{f_{i}\text{ continuous}}{=}
\end{equation*}%
\begin{equation*}
=\underset{n\rightarrow \infty }{\lim }f_{i}(f_{[\alpha ]_{n}}(x))=\underset{%
n\rightarrow \infty }{\lim }f_{i[\alpha ]_{n}}(x)=\underset{n\rightarrow
\infty }{\lim }f_{[i\alpha ]_{n+1}}(x)=a_{i\alpha }(x)\text{,}
\end{equation*}%
for all $x\in X$, $i\in I\cup J$\ and $\alpha \in \Lambda (I\cup J)$. $%
\square $

\bigskip

\textbf{Proposition 3.6.} \textit{Given an mIIFS} $\mathcal{S}%
=((X,d),(f_{i})_{i\in I\cup J})$, \textit{we have}%
\begin{equation*}
\underset{n\rightarrow \infty }{\lim }h^{\ast }(f_{[\alpha
]_{n}}(B),a_{\alpha }(B))=0\text{,}
\end{equation*}%
\textit{for all }$B\in P_{b}(X)$\textit{\ and }$\alpha \in \Lambda (I\cup J)$%
\textit{.}

\textit{Proof}. Let us consider $B\in P_{b}(X)$\textit{\ }and $\alpha \in
\Lambda (I\cup J)$ arbitrarily chosen, but fixed.

First let us note that 
\begin{equation*}
a_{\alpha }(B)\in P_{b}(X)
\end{equation*}%
and%
\begin{equation*}
f_{[\alpha ]_{n}}(B)\in P_{b}(X)\text{,}
\end{equation*}%
for every $n\in \mathbb{N}$, since $B\in P_{b}(X)$ and $a_{\alpha }$ and $%
f_{[\alpha ]_{n}}$ are Lipschitz.

Finally, we have%
\begin{equation*}
h^{\ast }(f_{[\alpha ]_{n}}(B),a_{\alpha }(B))\overset{\text{Proposition
2.2, iii)}}{\leq }
\end{equation*}%
\begin{equation*}
\leq \underset{x\in B}{\sup }\text{ }h^{\ast }(f_{[\alpha
]_{n}}(\{x\}),a_{\alpha }(\{x\}))\overset{\text{Proposition 2.2, i)}}{=}%
\underset{x\in B}{\sup }\text{ }d(f_{[\alpha ]_{n}}(x),a_{\alpha }(x))\text{,%
}
\end{equation*}%
for every $n\in \mathbb{N}$, so, using the squeeze theorem, via Corollary
3.2, we get the conclusion. $\square $

\bigskip

\textbf{Remark 3.7.} \textit{Given an mIIFS} $\mathcal{S}=((X,d),(f_{i})_{i%
\in I\cup J})$, \textit{it follows from Claim 2 from the proof of Theorem
3.7 from [1]} \textit{that} 
\begin{equation*}
\underset{n\rightarrow \infty }{\lim }\text{ }\underset{x\in B}{\sup }\text{ 
}h(F_{\mathcal{S}}^{[n]}(\{x\}),A_{x})=0\text{,}
\end{equation*}%
\textit{for every} $B\in P_{b,cl}(X)$\textit{.}

\bigskip

\textbf{Proposition 3.8.} \textit{Given an mIIFS} $\mathcal{S}%
=((X,d),(f_{i})_{i\in I\cup J})$, \textit{we have}%
\begin{equation*}
A_{B}=\overline{\underset{x\in B}{\cup }A_{x}}\text{,}
\end{equation*}%
\textit{for all }$B\in P_{b,cl}(X)$\textit{.}

\textit{Proof}. Let us consider $B\in P_{b,cl}(X)$ arbitrarily chosen, but
fixed.

In view of Proposition 2.8, we have%
\begin{equation*}
\underset{x\in B}{\cup }A_{x}\subseteq A_{B}\text{, }
\end{equation*}%
so $\underset{x\in B}{\cup }A_{x}\in P_{b}(X)$.

We have%
\begin{equation*}
h(A_{B},\overline{\underset{x\in B}{\cup }A_{x}})\leq h(A_{B},F_{\mathcal{S}%
}^{[n]}(B))+h^{\ast }(F_{\mathcal{S}}^{[n]}(B),\underset{x\in B}{\cup }%
A_{x})\leq
\end{equation*}%
\begin{equation*}
\overset{\text{Proposition 2.2, iii)}}{\leq }h(A_{B},F_{\mathcal{S}%
}^{[n]}(B))+\underset{x\in B}{\sup }\text{ }h(F_{\mathcal{S}%
}^{[n]}(\{x\}),A_{x})\text{,}
\end{equation*}%
for all $n\in \mathbb{N}$, so using the squeeze theorem, via Remark 3.7, we
get the conclusion. $\square $

\bigskip

\textbf{Proposition 3.9.} \textit{Let us consider an mIIFS }$\mathcal{S}%
=((X,d),(f_{i})_{i\in I\cup J})$, $\alpha =\alpha _{1}\alpha _{2}...\alpha
_{n}...\in \Lambda (I\cup J)$ \textit{and }$B\in P_{b,cl}(X)$\textit{\ such
that }$F_{\mathcal{S}}(B)\subseteq B$.

\textit{a) If }$\alpha \in \Sigma _{1}(I,J)$\textit{, i.e. }$n_{I}(\alpha
)=\infty $\textit{, then}%
\begin{equation*}
\underset{n\in \mathbb{N}}{\cap }\overline{f_{[\alpha ]_{n}}(B)}=Im%
a_{\alpha }
\end{equation*}%
\textit{and}%
\begin{equation*}
\underset{n\rightarrow \infty }{\lim }h(\overline{f_{[\alpha ]_{n}}(B)},%
Ima_{\alpha })=0\text{.}
\end{equation*}

\textit{b) If }$\alpha \in \Lambda (I\cup J)\setminus \Sigma _{1}(I,J)$%
\textit{, i.e. }$n_{I}(\alpha )<\infty $\textit{, then}%
\begin{equation*}
a_{\alpha }(B)=f_{\alpha _{1}\alpha _{2}...\alpha _{n^{\ast }}}(a_{\alpha
_{n^{\ast }+1}...\alpha _{m}...}(B))\text{,}
\end{equation*}%
\textit{and}%
\begin{equation*}
\overline{a_{\alpha _{n^{\ast }+1}...\alpha _{m}...}(B)}=\underset{p\in 
\mathbb{N}}{\cap }\overline{f_{\alpha _{n^{\ast }+1}...\alpha _{n^{\ast
}+p}}(B)}=\underset{p\rightarrow \infty }{\lim }\overline{f_{\alpha
_{n^{\ast }+1}...\alpha _{n^{\ast }+p}}(B)}\text{,}
\end{equation*}%
\textit{where} $n^{\ast }=\{%
\begin{array}{cc}
\max \{n\in \mathbb{N}\mid \alpha _{n}\in I\}\text{,} & \text{\textit{if} }%
n_{I}(\alpha )\neq 0 \\ 
0\text{,} & \text{\textit{otherwise}}%
\end{array}%
$,\textit{\ with the convention that }$f_{\alpha _{1}\alpha _{2}...\alpha
_{n^{\ast }}}=Id_{X}$\textit{\ if }$n^{\ast }=0$.

\textit{Proof}.

a) We can prove, via the mathematical induction method, that:

-%
\begin{equation*}
f_{[\alpha ]_{n+1}}(B)\subseteq f_{[\alpha ]_{n}}(B)\text{,}
\end{equation*}%
for every $n\in \mathbb{N}$;

-%
\begin{equation*}
diam(f_{[\alpha ]_{n}}(B))\leq a^{n_{I}([\alpha ]_{n})}diam(B)\text{,}
\end{equation*}%
for every $n\in \mathbb{N}$.

Consequently $\underset{n\in \mathbb{N}}{\cap }\overline{f_{[\alpha ]_{n}}(B)%
}$ consists of a single point.

Since $f_{[\alpha ]_{n+p}}(x)\in f_{[\alpha ]_{n}}(B)$ for every $x\in B$
and $n,p\in \mathbb{N}$, we infer that 
\begin{equation*}
a_{\alpha }(x)=\underset{p\rightarrow \infty }{\lim }f_{[\alpha
]_{n+p}}(x)\in \overline{f_{[\alpha ]_{n}}(B)}\text{,}
\end{equation*}%
for every $n\in \mathbb{N}$, so%
\begin{equation*}
a_{\alpha }(x)\in \underset{n\in \mathbb{N}}{\cap }\overline{f_{[\alpha
]_{n}}(B)}\text{.}
\end{equation*}

As, taking into account Proposition 3.4, b), $a_{\alpha }$ is constant, we
have 
\begin{equation*}
\underset{n\in \mathbb{N}}{\cap }\overline{f_{[\alpha ]_{n}}(B)}=Im%
a_{\alpha }\text{.}
\end{equation*}

Since 
\begin{equation*}
\underset{n\rightarrow \infty }{\lim }h(\overline{f_{[\alpha ]_{n}}(B)},%
\underset{n\in \mathbb{N}}{\cap }\overline{f_{[\alpha ]_{n}}(B)})\overset{%
\text{Proposition 2.5}}{=}0\text{,}
\end{equation*}%
we conclude that%
\begin{equation*}
\underset{n\rightarrow \infty }{\lim }h(\overline{f_{[\alpha ]_{n}}(B)},%
Ima_{\alpha })=0\text{.}
\end{equation*}

b) We have%
\begin{equation*}
a_{\alpha }(B)=\{a_{\alpha _{1}\alpha _{2}...\alpha _{n^{\ast }}\alpha
_{n^{\ast }+1}...\alpha _{m}...}(x)\mid x\in B\}\overset{\text{Proposition
3.5}}{=}
\end{equation*}%
\begin{equation*}
=\{f_{\alpha _{1}\alpha _{2}...\alpha _{n^{\ast }}}(a_{\alpha _{n^{\ast
}+1}...\alpha _{m}...}(x))\mid x\in B\}=f_{\alpha _{1}\alpha _{2}...\alpha
_{n^{\ast }}}(a_{\alpha _{n^{\ast }+1}...\alpha _{m}...}(B))\text{.}
\end{equation*}

As $\alpha _{k}\in J$ for every $k\in \mathbb{N}$, $k\geq n^{\ast }+1$, via
[16], we deduce that%
\begin{equation*}
\overline{a_{\alpha _{n^{\ast }+1}...\alpha _{m}...}(B)}=\underset{m\in 
\mathbb{N}}{\cap }\overline{f_{\alpha _{n^{\ast }+1}...\alpha _{n^{\ast
}+m}}(B)}=\underset{m\rightarrow \infty }{\lim }\overline{f_{\alpha
_{n^{\ast }+1}...\alpha _{n^{\ast }+m}}(B)}\text{. }\square 
\end{equation*}

\bigskip

\textbf{Lemma 3.10.} \textit{For every mIIFS }$\mathcal{S=}%
((X,d),(f_{i})_{i\in I\cup J})$\textit{, we have}%
\begin{equation*}
A_{a_{\alpha }(x)}\subseteq A_{x}\text{,}
\end{equation*}%
\textit{for all }$x\in X$\textit{\ and }$\alpha \in \Lambda (I\cup J)$%
\textit{.}

\textit{Proof}. Since $\underset{n\rightarrow \infty }{\lim }f_{[\alpha
]_{n}}(x)=a_{\mathcal{\alpha }}(x)$, we obtain%
\begin{equation*}
A_{a_{\alpha }(x)}\overset{\text{Lemma 2.13}}{=}\underset{n\rightarrow
\infty }{\lim }A_{f_{[\alpha ]_{n}}(x)}\overset{\text{Proposition 2.4}}{=}%
\{y\in X\mid \text{there exist }
\end{equation*}%
\begin{equation*}
\text{a strictly increasing sequence }(n_{k})_{k\in \mathbb{N}}\text{ of
natural numbers and}
\end{equation*}%
\begin{equation*}
y_{n_{k}}\in A_{f_{[\alpha ]_{n_{k}}}(x)}\overset{\text{Lemma 2.14}}{%
\subseteq }A_{x}\text{ for every }k\in \mathbb{N}\text{ such that }y=%
\underset{k\rightarrow \infty }{\lim }y_{n_{k}}\}\subseteq
\end{equation*}%
\begin{equation*}
\subseteq \overline{A_{x}}=A_{x}\text{,}
\end{equation*}%
for all $x\in X$ and $\alpha \in \Lambda (I\cup J)$.\textit{\ }$\square $

\bigskip

For an mIIFS $\mathcal{S=}((X,d),(f_{i})_{i\in I\cup J})$, $B\in P_{b,cl}(X)$
and $\alpha \in \Sigma (I,J)$, we introduce the function $\mathcal{A}%
_{\alpha }:X\rightarrow X$ as follows:

-%
\begin{equation*}
\mathcal{A}_{\alpha }=a_{\alpha }\text{,}
\end{equation*}%
if $\alpha \in \Sigma _{1}(I,J)$

-%
\begin{equation*}
\mathcal{A}_{\alpha }=f_{\beta _{0}}\circ a_{\gamma _{1}}\circ f_{\beta
_{1}}\circ a_{\gamma _{2}}\circ ...\circ f_{\beta _{n-1}}\circ a_{\gamma
_{n}}\circ f_{\beta _{n}}\text{,}
\end{equation*}%
if $\alpha =\beta _{0}\gamma _{1}\beta _{1}...\gamma _{n}\beta _{n}\in
\Sigma _{0}(I,J)$, where $n\in \mathbb{N}$, $\beta _{0}\in \{\lambda \}\cup
\Lambda _{1}(I\cup J)$, $\beta _{n}\in \{\lambda \}\cup \Lambda _{2}(I\cup J)
$, $\gamma _{k}\in \Lambda (J)$ for every $k\in \{1,...,n\}$ and $\beta
_{k}\in \Lambda _{3}(I\cup J)$ for every $k\in \{1,...,n-1\}$ if $n\geq 2$.

Note that $\mathcal{A}_{\alpha }$ is well defined.

We also consider%
\begin{equation*}
\{\mathcal{A}_{\alpha }(x)\mid \alpha \in \Sigma (I,J)\text{ and }x\in B\}%
\overset{not}{=}L_{B}\text{.}
\end{equation*}

Note that, in view of the next Lemma, $L_{B}\in P_{b}(X)$.

\bigskip

\textbf{Lemma 3.11.} \textit{For every mIIFS }$\mathcal{S=}%
((X,d),(f_{i})_{i\in I\cup J})$ \textit{and} $B\in P_{b,cl}(X)$\textit{, we
have}%
\begin{equation*}
\overline{L_{B}}\subseteq A_{B}\text{.}
\end{equation*}

\textit{Proof}. One the one hand, we have%
\begin{equation}
\mathcal{A}_{\alpha }(x)\overset{\text{Proposition 3.3}}{\in }A_{B}\text{,} 
\tag{1}
\end{equation}%
for all $x\in B$ and $\alpha \in \Sigma _{1}(I,J)$.

On the other hand, for all $x\in B$ and $\alpha =\beta _{0}\gamma _{1}\beta
_{1}...\gamma _{n}\beta _{n}\in \Sigma _{0}(I,J)$, where $n\in \mathbb{N}$, $%
\beta _{0}\in \{\lambda \}\cup \Lambda _{1}(I\cup J)$, $\beta _{n}\in
\{\lambda \}\cup \Lambda _{2}(I\cup J)$, $\gamma _{k}\in \Lambda (J)$ for
every $k\in \{1,...,n\}$ and $\beta _{k}\in \Lambda _{3}(I\cup J)$ for every 
$k\in \{1,...,n-1\}$ if $n\geq 2$, we have 
\begin{equation*}
a_{\gamma _{_{n}}}(f_{\beta _{n}}(x))\overset{\text{Proposition 3.3}}{\in }%
A_{f_{\beta _{n}}(x)}\overset{\text{Lemma 2.14}}{\subseteq }A_{x}\text{,}
\end{equation*}%
so 
\begin{equation*}
a_{\gamma _{_{n-1}}}(f_{\beta _{n-1}}(a_{\gamma _{_{n}}}(f_{\beta _{n}}(x))))%
\overset{\text{Proposition 3.3}}{\in }A_{f_{\beta _{n-1}}(a_{\gamma
_{_{n}}}(f_{\beta _{n}}(x)))}\subseteq 
\end{equation*}%
\begin{equation*}
\overset{\text{Lemma 2.14}}{\subseteq }A_{a_{\gamma _{_{n}}}(f_{\beta
_{n}}(x))}\overset{\text{Lemma 3.10}}{\subseteq }A_{f_{\beta _{n}}(x)}%
\overset{\text{Lemma 2.14}}{\subseteq }A_{x}\text{,}
\end{equation*}%
and, inductively, we obtain%
\begin{equation*}
(f_{\beta _{0}}\circ a_{\gamma _{1}}\circ f_{\beta _{1}}\circ a_{\gamma
_{2}}\circ ...\circ f_{\beta _{n-1}}\circ a_{\gamma _{n}}\circ f_{\beta
_{n}})(x)\overset{\text{Proposition 3.5}}{=}
\end{equation*}%
\begin{equation*}
=(a_{\beta _{0}\gamma _{1}}\circ f_{\beta _{1}}\circ a_{\gamma _{2}}\circ
...\circ f_{\beta _{n-1}}\circ a_{\gamma _{n}}\circ f_{\beta _{n}})(x)\in
A_{x}\text{.}
\end{equation*}

Hence%
\begin{equation}
\mathcal{A}_{\alpha }(x)\in A_{x}\overset{\text{Proposition 3.8}}{\subseteq }%
A_{B}\text{,}  \tag{2}
\end{equation}%
for all $x\in B$ and $\alpha \in \Sigma _{0}(I,J)$.

Therefore, via $(1)$ and $(2)$, we get $L_{B}\subseteq A_{B}$, and
consequently $\overline{L_{B}}\subseteq A_{B}$. $\square $

\bigskip

The next result provides an alternative description of the fixed points of
the fractal operator associated with an mIIFS.

\bigskip

\textbf{Theorem 3.12.} \textit{For every mIIFS }$\mathcal{S=}%
((X,d),(f_{i})_{i\in I\cup J})$ \textit{and} $B\in P_{b,cl}(X)$\textit{, we
have}%
\begin{equation*}
\overline{L_{B}}=A_{B}\text{,}
\end{equation*}%
\textit{i.e.}%
\begin{equation*}
\overline{\{\mathcal{A}_{\alpha }(x)\mid \alpha \in \Sigma (I,J)\text{ 
\textit{and} }x\in B\}}=A_{B}\text{.}
\end{equation*}

\textit{Proof}. Let $B\in P_{b,cl}(X)$ be arbitrarily chosen, but fixed.

We will prove that%
\begin{equation}
D(A_{B},\overline{L_{B}})=0\text{.}  \tag{*}
\end{equation}

In view of Remark 3.7 and of the following inequality 
\begin{equation*}
D(A_{B},\overline{L_{B}})\overset{\text{Proposition 3.8 }}{=}D(\overline{%
\underset{x\in B}{\cup }A_{x}},\overline{L_{B}})\overset{\text{Remark 2.1, i)%
}}{=}
\end{equation*}%
\begin{equation*}
=D(\underset{x\in B}{\cup }A_{x},\overline{L_{B}})\overset{\text{Remark 2.1,
ii)}}{\leq }\underset{x\in B}{\sup }\text{ }D(A_{x},\overline{L_{B}})\leq 
\end{equation*}%
\begin{equation*}
\leq \underset{x\in B}{\sup }\text{ }D(A_{x},F_{\mathcal{S}}^{[n]}(\{x\}))+%
\underset{x\in B}{\sup }\text{ }D(F_{\mathcal{S}}^{[n]}(\{x\}),\overline{%
L_{B}})\leq 
\end{equation*}%
\begin{equation*}
\leq \underset{x\in B}{\sup }\text{ }h(A_{x},F_{\mathcal{S}}^{[n]}(\{x\}))+%
\underset{x\in B}{\sup }\text{ }D(F_{\mathcal{S}}^{[n]}(\{x\}),\overline{%
L_{B}})\text{,}
\end{equation*}%
which is valid for all $n\in \mathbb{N}$, it suffices to prove that%
\begin{equation}
\underset{n\rightarrow \infty }{\lim }\text{ }\underset{x\in B}{\sup }\text{ 
}D(F_{\mathcal{S}}^{[n]}(\{x\}),\overline{L_{B}})=0\text{.}  \tag{1}
\end{equation}

In order to justify $(1)$, let us consider $\varepsilon >0$ fixed, but
arbitrarily chosen.

Then there exists $n_{1}\in \mathbb{N}$ such that 
\begin{equation}
\frac{a^{[\sqrt{n}]-2}}{1-a}N_{B}<\frac{\varepsilon }{2}\text{,}  \tag{2}
\end{equation}%
for every $n\in \mathbb{N}$, $n\geq n_{1}$, where $a\in \lbrack 0,1)$ is the
constant associated with $S$ via Definition 2.11

Now let us consider $n\in \mathbb{N}$, $n\geq n_{1}$, $x\in B$ and $z\in F_{%
\mathcal{S}}^{[n]}(\{x\})$.

Taking into account Proposition 2.8, there exists $\alpha =\alpha
_{1}...\alpha _{n}\in \Lambda ^{\ast }(I\cup J)$ such that%
\begin{equation}
d(z,f_{\alpha }(x))<\frac{\varepsilon }{2}\text{.}  \tag{3}
\end{equation}

The following two cases occur:

i) $n_{I}(\alpha )\geq \lbrack \sqrt{n}]$

ii) $n_{I}(\alpha )<[\sqrt{n}]$.

In the first case, we consider $\gamma \in \Lambda (I)$ and we have%
\begin{equation}
d(f_{[\alpha \gamma ]_{\left\vert \alpha \right\vert +k}}(x),f_{[\alpha
\gamma ]_{\left\vert \alpha \right\vert +k+1}}(x))\overset{\text{Definition
2.11, b 1)}}{\leq }a^{n_{I}(\alpha )+k}N_{B}\text{,}  \tag{4}
\end{equation}%
for every $k\in \mathbb{N}$, so%
\begin{equation*}
d(f_{\alpha }(x),\mathcal{A}_{\alpha \gamma }(x))\leq d(f_{\alpha
}(x),f_{[\alpha \gamma ]_{\left\vert \alpha \right\vert
+l+1}}(x))+d(f_{[\alpha \gamma ]_{\left\vert \alpha \right\vert +l+1}}(x),%
\mathcal{A}_{\alpha \gamma }(x))\leq
\end{equation*}%
\begin{equation*}
\leq \overset{l}{\underset{k=0}{\sum }}d(f_{[\alpha \gamma ]_{\left\vert
\alpha \right\vert +k}}(x),f_{[\alpha \gamma ]_{\left\vert \alpha
\right\vert +k+1}}(x))+d(f_{[\alpha \gamma ]_{\left\vert \alpha \right\vert
+l+1}}(x),\mathcal{A}_{\alpha \gamma }(x))\overset{\text{(4)}}{\leq }
\end{equation*}%
\begin{equation*}
\leq a^{n_{I}(\alpha )}N_{B}\overset{l}{\underset{k=0}{\sum }}%
a^{k}+d(f_{[\alpha \gamma ]_{\left\vert \alpha \right\vert +l+1}}(x),%
\mathcal{A}_{\alpha \gamma }(x))\leq
\end{equation*}%
\begin{equation*}
\leq \frac{a^{n_{I}(\alpha )}}{1-a}N_{B}+d(f_{[\alpha \gamma ]_{\left\vert
\alpha \right\vert +l+1}}(x),\mathcal{A}_{\alpha \gamma }(x))\text{,}
\end{equation*}%
for every $l\in \mathbb{N}$. By passing to the limit, as $l$ goes to $\infty 
$, in the previous inequality, we get%
\begin{equation}
d(f_{\alpha }(x),\mathcal{A}_{\alpha \gamma }(x))\leq \frac{a^{n_{I}(\alpha
)}}{1-a}N_{B}\leq \frac{a^{[\sqrt{n}]-2}}{1-a}N_{B}\overset{\text{(2)}}{<}%
\frac{\varepsilon }{2}\text{.}  \tag{5}
\end{equation}%
Hence we have%
\begin{equation*}
d(z,\overline{L_{B}})\overset{\mathcal{A}_{\alpha \gamma }(x)\in L_{B}}{\leq 
}d(z,\mathcal{A}_{\alpha \gamma }(x))\leq
\end{equation*}%
\begin{equation}
\leq d(z,f_{\alpha }(x))+d(f_{\alpha }(x),\mathcal{A}_{\alpha \gamma }(x))%
\overset{\text{(3) \& (5)}}{\leq }\frac{\varepsilon }{2}+\frac{\varepsilon }{%
2}=\varepsilon \text{.}  \tag{6}
\end{equation}

In the second case, there exist $k\in \mathbb{N}$, $\beta _{1},...,\beta
_{k+1}\in \Lambda ^{\ast }(I)$ and $\gamma _{j}\in \Lambda ^{\ast
}(J)\smallsetminus \{\lambda \}$, $j\in \{1,...,k\}$ such that 
\begin{equation*}
\alpha =\beta _{1}\gamma _{1}\beta _{2}\gamma _{2}...\beta _{k}\gamma
_{k}\beta _{k+1}\text{.}
\end{equation*}%
On the one hand, $\left\vert \beta _{1}\right\vert +\left\vert \beta
_{2}\right\vert +...+\left\vert \beta _{k+1}\right\vert =n_{I}(\alpha )<[%
\sqrt{n}]$ and $\left\vert \beta _{1}\right\vert +\left\vert \beta
_{2}\right\vert +...+\left\vert \beta _{k+1}\right\vert +\left\vert \gamma
_{1}\right\vert +\left\vert \gamma _{2}\right\vert +...+\left\vert \gamma
_{k}\right\vert =n$, so 
\begin{equation}
\left\vert \gamma _{1}\right\vert +\left\vert \gamma _{2}\right\vert
+...+\left\vert \gamma _{k}\right\vert >n-[\sqrt{n}]\text{.}  \tag{7}
\end{equation}%
On the other hand, as $k-1\leq n_{I}(\alpha )<[\sqrt{n}]$, we have%
\begin{equation}
k\leq \lbrack \sqrt{n}]+1\text{.}  \tag{8}
\end{equation}%
Based on $(7)$ and $(8)$, via the generalized pigeonhole principle, there
exists $j\in \{1,...,k\}$ such that%
\begin{equation*}
\left\vert \gamma _{j}\right\vert \geq \lbrack \frac{\left\vert \gamma
_{1}\right\vert +\left\vert \gamma _{2}\right\vert +...+\left\vert \gamma
_{k}\right\vert -1}{k}]+1\geq \frac{\left\vert \gamma _{1}\right\vert
+\left\vert \gamma _{2}\right\vert +...+\left\vert \gamma _{k}\right\vert -1%
}{k}\overset{\text{(7) \& (8)}}{\geq }
\end{equation*}%
\begin{equation}
\geq \frac{n-[\sqrt{n}]-1}{[\sqrt{n}]+1}\geq \frac{([\sqrt{n}])^{2}-[\sqrt{n}%
]-1}{[\sqrt{n}]+1}=[\sqrt{n}]-2+\frac{1}{[\sqrt{n}]+1}\geq \lbrack \sqrt{n}%
]-2\text{.}  \tag{9}
\end{equation}%
In the sequel, we shall use the following notation:%
\begin{equation*}
(f_{\beta _{j+1}}\circ f_{\gamma _{j+1}}\circ ...\circ f_{\beta _{k}}\circ
f_{\gamma _{k}}\circ f_{\beta _{k+1}})(x)\overset{not}{=}y\text{.}
\end{equation*}%
For a fixed $\gamma =\gamma ^{1}\gamma ^{2}...\gamma ^{n}...\in \Lambda (J)$%
, let us consider%
\begin{equation*}
\theta =\beta _{1}\gamma _{1}\beta _{2}\gamma _{2}...\beta _{j}\gamma
_{j}\gamma \beta _{j+1}\gamma _{j+1}...\beta _{k}\gamma _{k}\beta _{k+1}\in
\Sigma _{0}(I,J)
\end{equation*}%
and%
\begin{equation*}
\theta _{l}=\beta _{1}\gamma _{1}\beta _{2}\gamma _{2}...\beta _{j}\gamma
_{j}[\gamma ]_{l}\beta _{j+1}\gamma _{j+1}...\beta _{k}\gamma _{k}\beta
_{k+1}\text{,}
\end{equation*}%
where $l\in \{0\}\cup \mathbb{N}$. Then%
\begin{equation*}
d(f_{\alpha }(x),\mathcal{A}_{\theta }(x))\leq d(f_{\theta
_{0}}(x),f_{\theta _{l}}(x))+d(f_{\theta _{l}}(x),\mathcal{A}_{\theta
}(x))\leq 
\end{equation*}%
\begin{equation*}
\leq \overset{l-1}{\underset{s=0}{\sum }}d(f_{\theta _{s}}(x),f_{\theta
_{s+1}}(x))+d(f_{\theta _{l}}(x),\mathcal{A}_{\theta }(x))\overset{\text{(9)}%
}{\leq }
\end{equation*}%
\begin{equation*}
\leq \overset{l-1}{\underset{s=0}{\sum }}a^{[\sqrt{n}]-2+s}d(y,f_{\gamma
^{s+1}}(y))+d(f_{\theta _{l}}(x),\mathcal{A}_{\theta }(x))\leq 
\end{equation*}%
\begin{equation*}
\leq \overset{l-1}{\underset{s=0}{\sum }}a^{[\sqrt{n}]-2+s}N_{B}+
\end{equation*}%
\begin{equation*}
+d((f_{\beta _{1}}\circ f_{\gamma _{1}}\circ ...\circ f_{\gamma _{j}}\circ
f_{[\gamma ]_{l}})(y),(f_{\beta _{1}}\circ f_{\gamma _{1}}\circ ...\circ
f_{\gamma _{j}}\circ a_{\gamma })(y))\leq 
\end{equation*}%
\begin{equation}
\leq N_{B}\overset{l-1}{\underset{s=0}{\sum }}a^{[\sqrt{n}%
]-2+s}+d(f_{[\gamma ]_{l}}(y),a_{\gamma }(y))\leq \frac{a^{[\sqrt{n}]-2}}{1-a%
}N_{B}+d(f_{[\gamma ]_{l}}(y),a_{\gamma }(y))\text{,}  \tag{10}
\end{equation}%
for every $l\in \{0\}\cup \mathbb{N}$, so, by passing to limit, as $l$ goes
to $\infty $, in $(10)$, we infer that%
\begin{equation}
d(f_{\alpha }(x),\mathcal{A}_{\theta }(x))\leq \frac{a^{[\sqrt{n}]-2}}{1-a}%
N_{B}\overset{\text{(2)}}{<}\frac{\varepsilon }{2}\text{.}  \tag{11}
\end{equation}%
As in the first case, we have%
\begin{equation*}
d(z,\overline{L_{B}})\overset{\mathcal{A}_{\theta }(x)\in L_{B}}{\leq }d(z,%
\mathcal{A}_{\theta }(x))\leq 
\end{equation*}%
\begin{equation}
\leq d(z,f_{\alpha }(x))+d(f_{\alpha }(x),\mathcal{A}_{\theta }(x))\overset{%
\text{(3) \& (11)}}{\leq }\frac{\varepsilon }{2}+\frac{\varepsilon }{2}%
=\varepsilon \text{.}  \tag{12}
\end{equation}

Via $(6)$ and $(12)$, we conclude that%
\begin{equation*}
D(F_{S}^{[n]}(x),\overline{L_{B}})=\underset{z\in F_{S}^{[n]}(x)}{\sup }d(z,%
\overline{L_{B}})\leq \varepsilon \text{,}
\end{equation*}%
for every $n\in \mathbb{N}$, $n\geq n_{1}$ and every $x\in B$, i.e.%
\begin{equation*}
\underset{n\rightarrow \infty }{\lim }\text{ }\underset{x\in B}{\sup }\text{ 
}D(F_{S}^{[n]}(x),\overline{L_{B}})=0\text{,}
\end{equation*}%
so $(\ast )$ is proved and therefore%
\begin{equation*}
A_{B}\subseteq \overline{L_{B}}\text{.}
\end{equation*}

The above inclusion, together with Lemma 3.11, completes the proof. $\square 
$

\bigskip

\textbf{4. The canonical projection associated with an mIIFS}

\bigskip

Finally let us introduce the concept of canonical projection associated with
an mIIFS and rewrite part of the obtained results.

\bigskip

\textbf{Definition 4.1.} \textit{For an} \textit{mIIFS} $\mathcal{S=}%
((X,d),(f_{i})_{i\in I\cup J})$\textit{, the function }$\pi :\Sigma
(I,J)\times X\rightarrow X$\textit{, given by}%
\begin{equation*}
\pi (\alpha ,x)=\mathcal{A}_{\alpha }(x)\text{,}
\end{equation*}%
\textit{for every }$(\alpha ,x)\in \Sigma (I,J)\times X$\textit{, is called} 
\textit{the canonical projection associated with }$\mathcal{S}$\textit{.}

\bigskip

Proposition 3.4, b) takes the following form (which is a counterpart of
Theorem 2.9, b)):

\bigskip

\textbf{Proposition 4.2.} \textit{For each} \textit{mIIFS} $\mathcal{S=}%
((X,d),(f_{i})_{i\in I\cup J})$,\textit{\ the set }$\pi (\alpha ,X)$\textit{%
\ has just one element for all }$\alpha \in \Sigma _{1}(I,J)$\textit{.}

\bigskip

The following result is the counterpart of Theorem 2.9, c), i).

\bigskip

\textbf{Proposition 4.3. }\textit{For each} \textit{mIIFS} $\mathcal{S=}%
((X,d),(f_{i})_{i\in I\cup J})$,\textit{\ we have}%
\begin{equation*}
\pi \circ \mathbf{\tau }_{i}=f_{i}\circ \pi \text{,}
\end{equation*}%
\textit{for every }$i\in I\cup J$\textit{, where }$\mathbf{\tau }_{i}:\Sigma
(I,J)\times X\rightarrow \Sigma (I,J)\times X$ \textit{is} \textit{defined by%
}%
\begin{equation*}
\mathbf{\tau }_{i}(\alpha ,x)=(i\alpha ,x)\text{,}
\end{equation*}%
\textit{for every }$(\alpha ,x)\in \Sigma (I,J)\times X$\textit{, i.e. the
following diagram is commutative}%
\begin{equation*}
\begin{array}{ccc}
\Sigma (I,J)\times X & \overset{\mathbf{\tau }_{i}}{\rightarrow } & \Sigma
(I,J)\times X \\ 
\pi \downarrow &  & \downarrow \pi \\ 
X & \overset{f_{i}}{\rightarrow } & X%
\end{array}%
\text{.}
\end{equation*}

\textit{Proof}. For $\alpha \in \Sigma _{1}(I,J)$ the conclusion results
from Proposition 3.5.

For $\alpha \in \Sigma _{0}(I,J)$, as $i\alpha \in \Sigma _{0}(I,J)$, the
conclusion takes the following obvious form:%
\begin{equation*}
(f_{i\beta _{0}}\circ a_{\gamma _{1}}\circ f_{\beta _{1}}\circ a_{\gamma
_{2}}\circ ...\circ f_{\beta _{n-1}}\circ a_{\gamma _{n}}\circ f_{\beta
_{n}})(x)=
\end{equation*}%
\begin{equation*}
=f_{i}((f_{\beta _{0}}\circ a_{\gamma _{1}}\circ f_{\beta _{1}}\circ
a_{\gamma _{2}}\circ ...\circ f_{\beta _{n-1}}\circ a_{\gamma _{n}}\circ
f_{\beta _{n}})(x))\text{,}
\end{equation*}%
where $\alpha =\beta _{0}\gamma _{1}\beta _{1}...\gamma _{n}\beta _{n}\in
\Sigma _{0}(I,J)$, with $n\in \mathbb{N}$, $\beta _{0}\in \{\lambda \}\cup
\Lambda _{1}(I\cup J)$, $\beta _{n}\in \{\lambda \}\cup \Lambda _{2}(I\cup J)
$, $\gamma _{k}\in \Lambda (J)$ for every $k\in \{1,...,n\}$ and $\beta
_{k}\in \Lambda _{3}(I\cup J)$ for every $k\in \{1,...,n-1\}$ if $n\geq 2$. $%
\square $

\bigskip

Theorem 3.12 can be rewritten in the following form, which is the
counterpart of Theorem 2.9, c) ii):

\bigskip

\textbf{Theorem 4.4. }\textit{Let }$\mathcal{S=}((X,d),(f_{i})_{i\in I\cup
J})$ \textit{be a} \textit{mIIFS} \textit{and} $B\in P_{b,cl}(X)$\textit{.
Then}%
\begin{equation*}
\overline{\pi (\Sigma (I,J)\times B)}=A_{B}\text{.}
\end{equation*}

\bigskip

\textbf{5. Visual aspects concerning the functions }$A_{\alpha }$

\bigskip

Let us consider the mIIFS $\mathcal{S}=((\mathbb{R}^{2},\left\Vert
.\right\Vert _{2}),(f_{i})_{i\in I\cup J})$, where $I=\{1,2,3\}$, $J=\{4\}$
and $f_{i}:\mathbb{R}^{2}\rightarrow \mathbb{R}^{2}$ are given by%
\begin{equation*}
f_{1}(x,y)=(\frac{x}{2},\frac{y}{2})\text{,}
\end{equation*}%
\begin{equation*}
f_{2}(x,y)=(\frac{x+1}{2},\frac{y}{2})\text{,}
\end{equation*}%
\begin{equation*}
f_{3}(x,y)=(\frac{2x+1}{4},\frac{2y+\sqrt{3}}{4})
\end{equation*}%
and%
\begin{equation*}
f_{4}(x,y)=(x,\frac{1}{5}y)\text{,}
\end{equation*}%
for every $(x,y)\in \mathbb{R}^{2}$.

\textbf{A}. We start by providing a visualization of the convergence of the
sequence defining $a_{\alpha }(u)$ for a randomly generated $\alpha \in
\Sigma _{1}(I,J)$ and for $u\in \{(0,0),(1,1),(20,20)\}$. More precisely, we
represent the first 50000 terms of the corresponding sequences, the limit
being marked by *.

\begin{figure}[H]
    \centering
\subfloat[\(a_\alpha((0,0))\)]{
\includegraphics[width=0.32\textwidth]{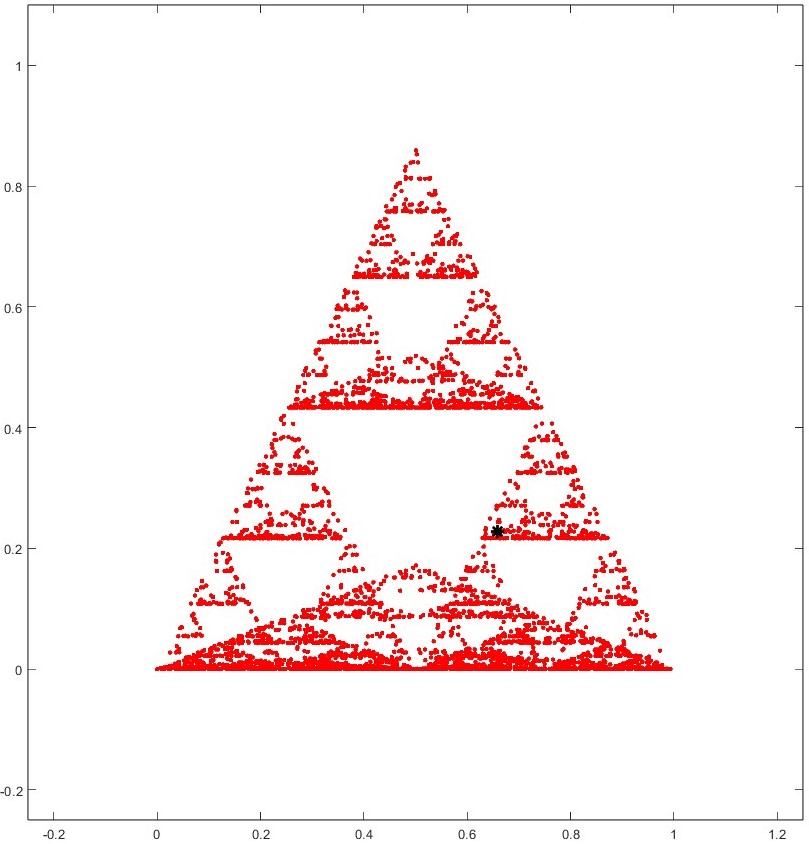}}
\subfloat[\(a_\alpha((1,1))\)]{
\includegraphics[width=0.33\textwidth]{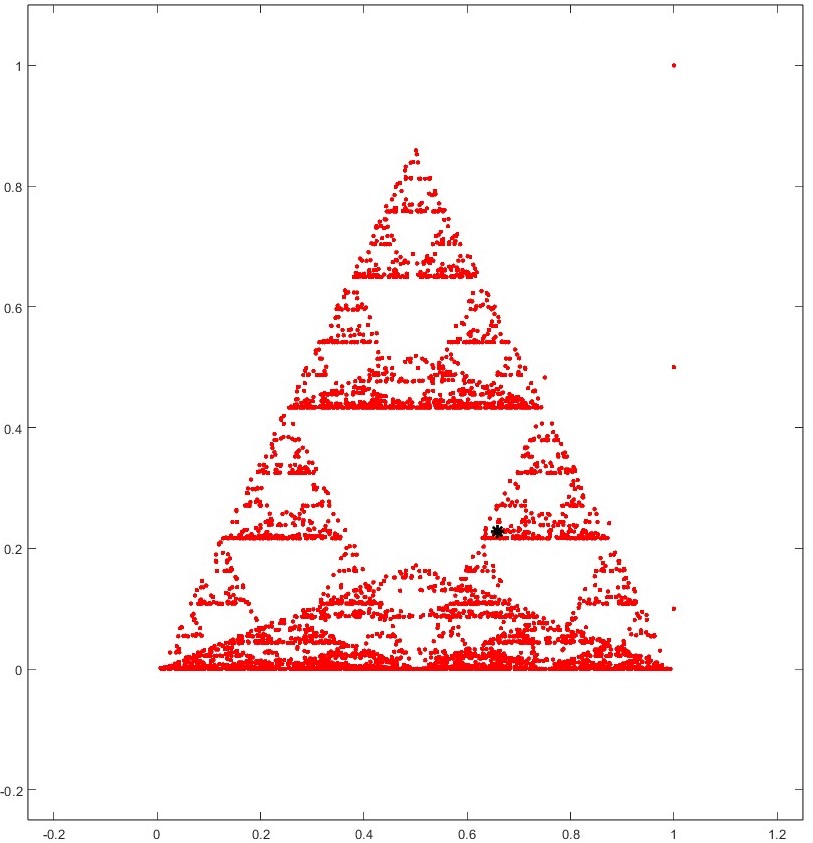}}
\subfloat[\(a_\alpha((20,20))\)]{
\includegraphics[width=0.33\textwidth]{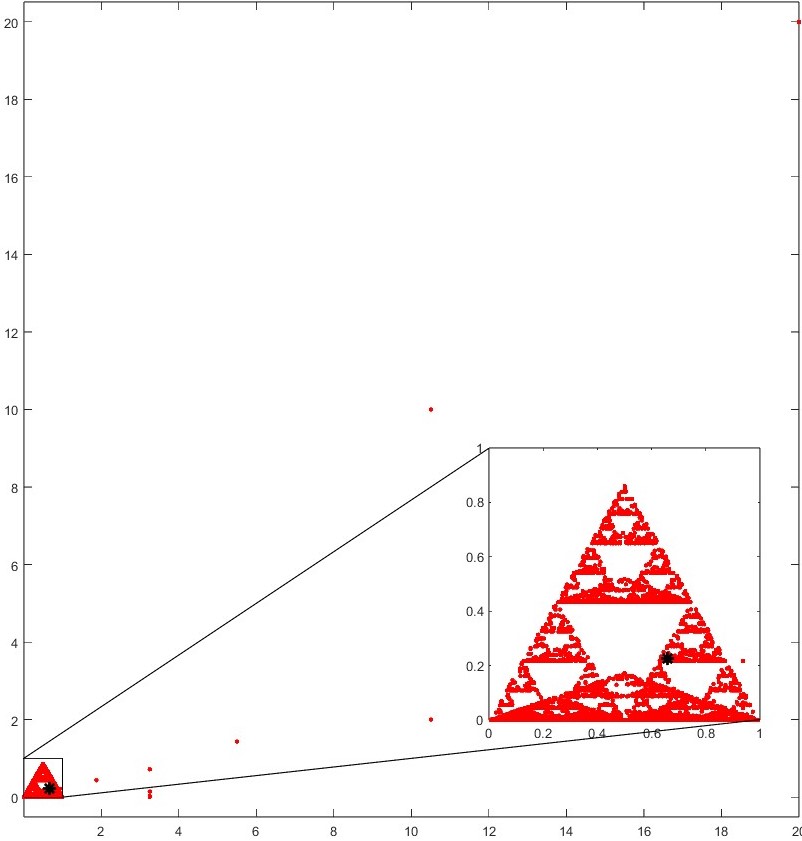}}
\caption{}
\end{figure}

Let us remark that the above presented visualizations are in accordance with
Proposition 3.4, b), as all the three pictures indicate a convergence to the
same limit.

\textbf{B}. Next we come up with a visualization of $\mathcal{A}_{\alpha
}(u) $ for $u\in \{(0,0),(1,1)\}$ and for $\alpha \in \{\omega ,321\omega
,\omega 123,12\omega 31\}\subseteq \Sigma _{1}(I,J)$, where $\omega $ is the
word having all letters equal to $4$.

\begin{figure}[H]
    \centering
\includegraphics[width=0.70\textwidth]{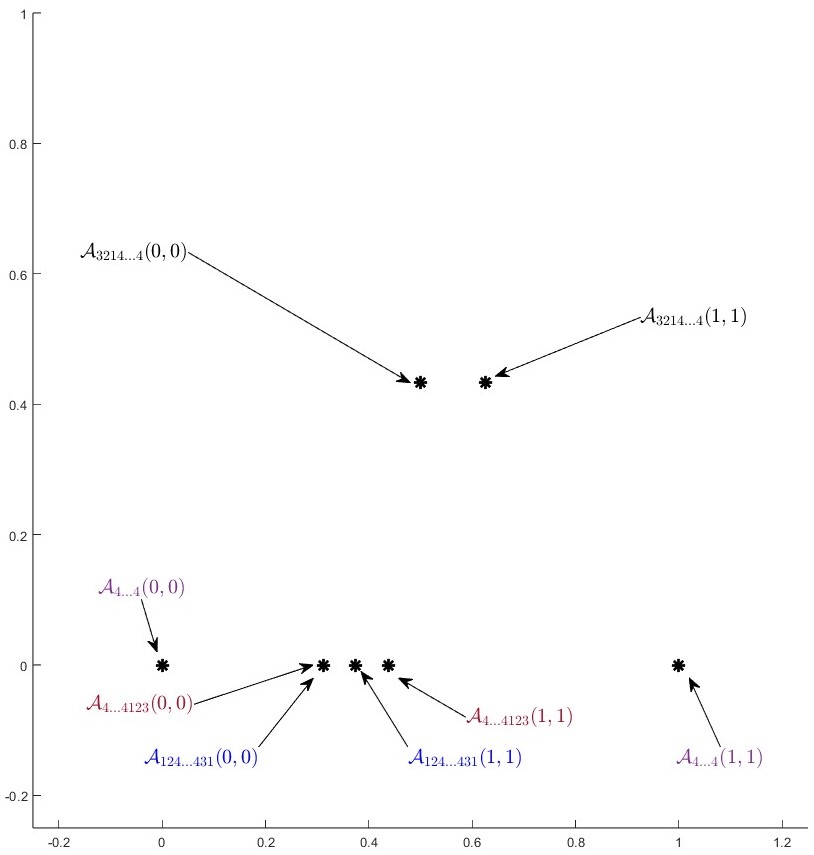}
\caption{}
\end{figure}
Let us remark that Figure 2 gives evidence that $\mathcal{A}_{\alpha }(u)$
depends on $u$.

\textbf{C}. Let us consider the set $B=F_{\mathcal{S}}^{[20]}(\{(0,0)\})$, $%
\alpha =1434121212...\in \Sigma _{1}(I,J)$ and $\beta =412\omega \in \Sigma
_{0}(I,J)$.

Figures 3(a), 3(b), 3(c), 3(d) and 3(e) contain the visual representations of $f_{[\alpha
]_{n}}(B)$ for $n=1$, $n=2$, $n=3$, $n=6$ and $n=9$, respectively. The fact
that $f_{[\alpha ]_{9}}(B)$ is very close with $Ima_{\alpha }$
(which is marked with *) endorses Theorem 3.9, a).
\begin{figure}[H]
    \centering
\subfloat[]{
\includegraphics[width=0.32\textwidth]{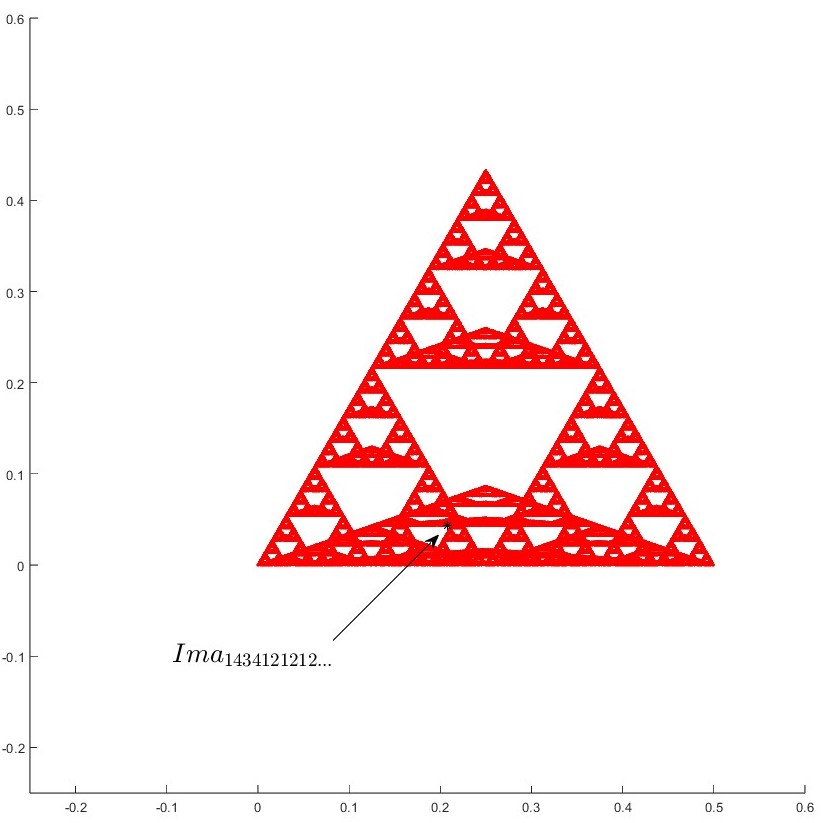}}
\subfloat[]{
\includegraphics[width=0.33\textwidth]{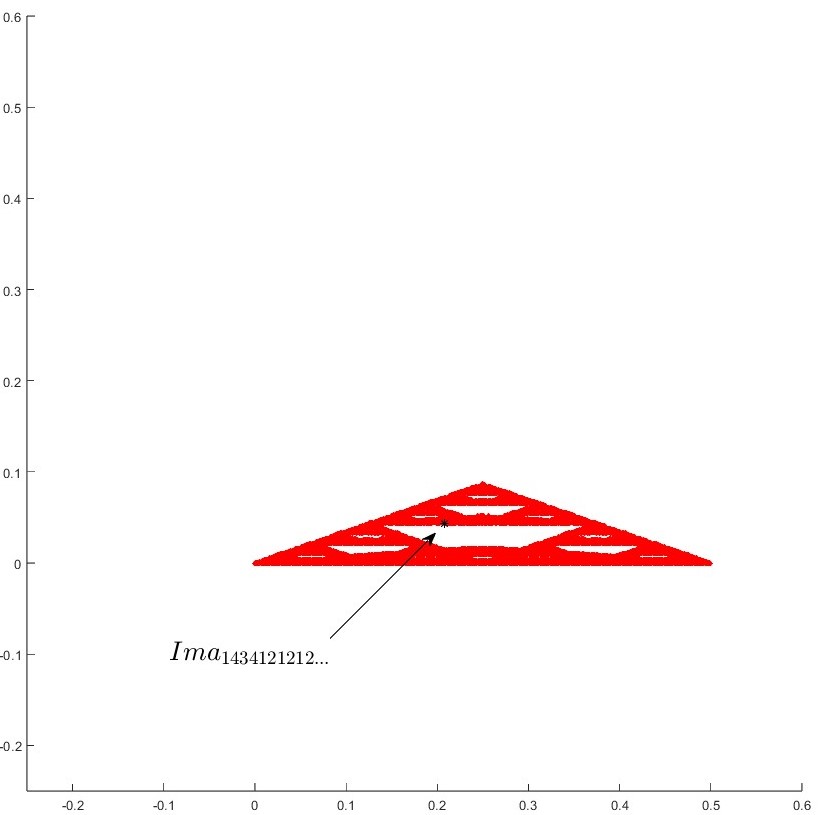}}
\subfloat[]{
\includegraphics[width=0.33\textwidth]{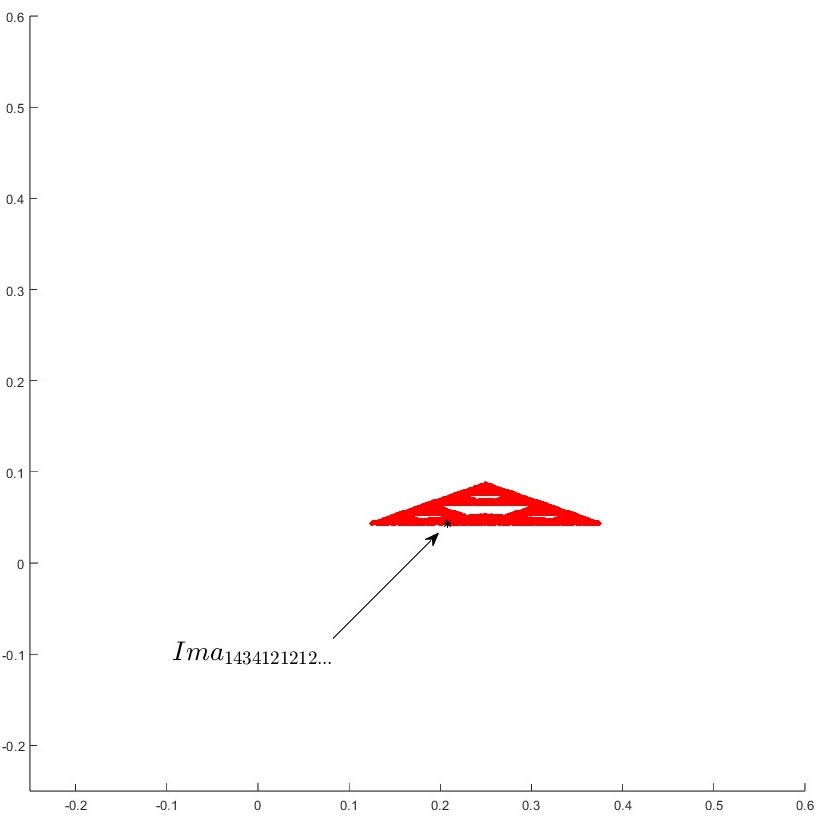}}\\
\subfloat[]{
\includegraphics[width=0.33\textwidth]{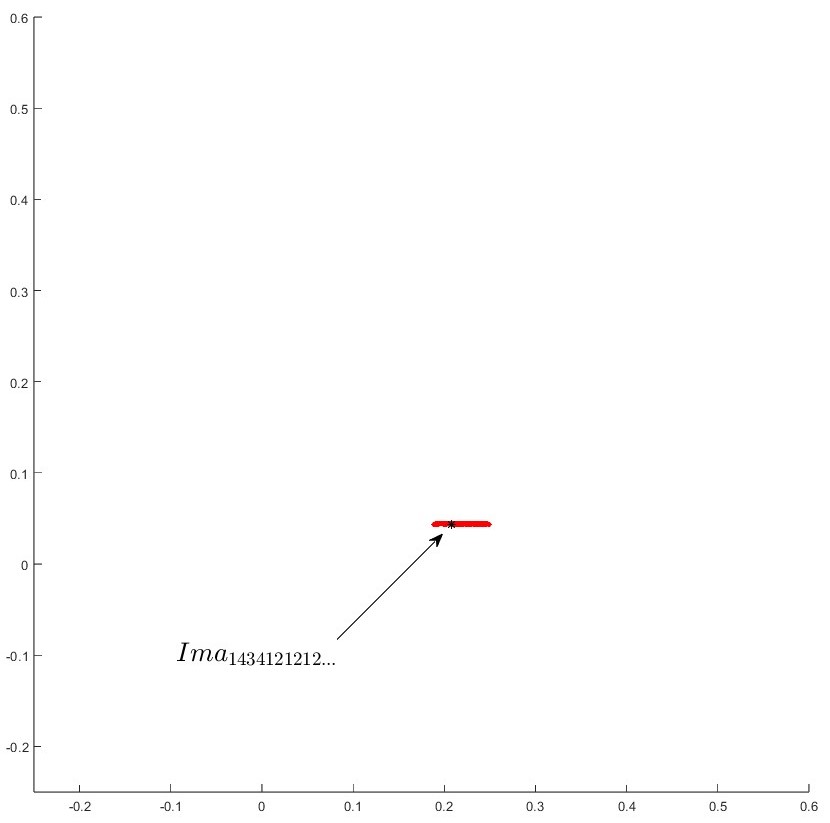}}
\subfloat[]{
\includegraphics[width=0.33\textwidth]{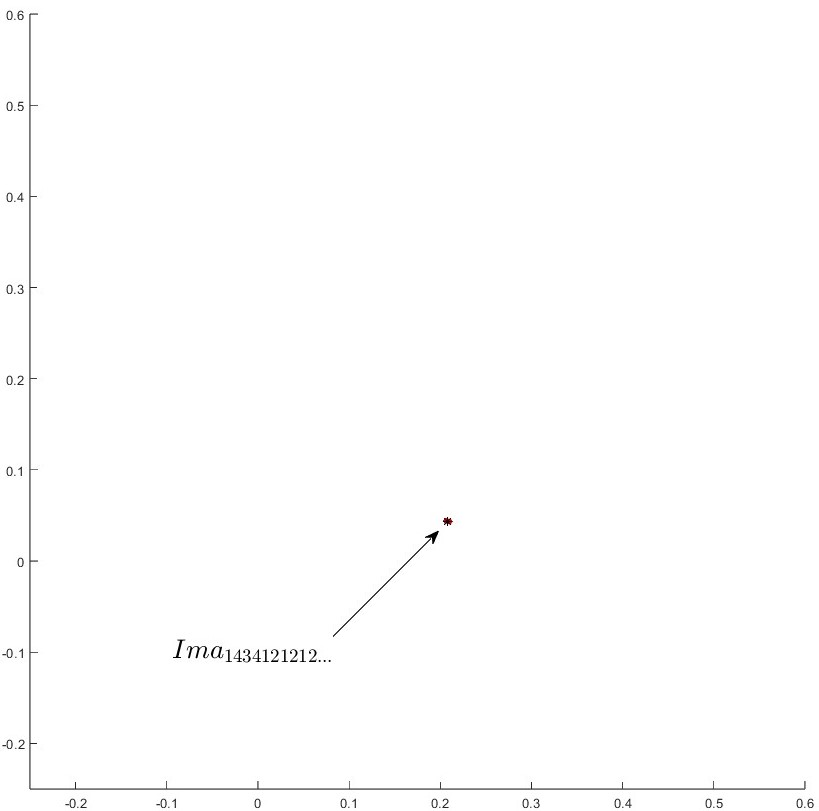}}
\caption{}
\end{figure}

Figures 4(a), 4(b), 4(c), 4(d) and 4(e) contain the visual representations of $%
f_{[\beta ]_{n}}(B)$ for $n=1$, $n=2$, $n=3$, $n=4$ and $n=6$, respectively
and they illustrate Theorem 3.9, b) (the black horizontal segment
representing $a_{\beta }(B)$).

\begin{figure}[H]
    \centering
\subfloat[]{
\includegraphics[width=0.32\textwidth]{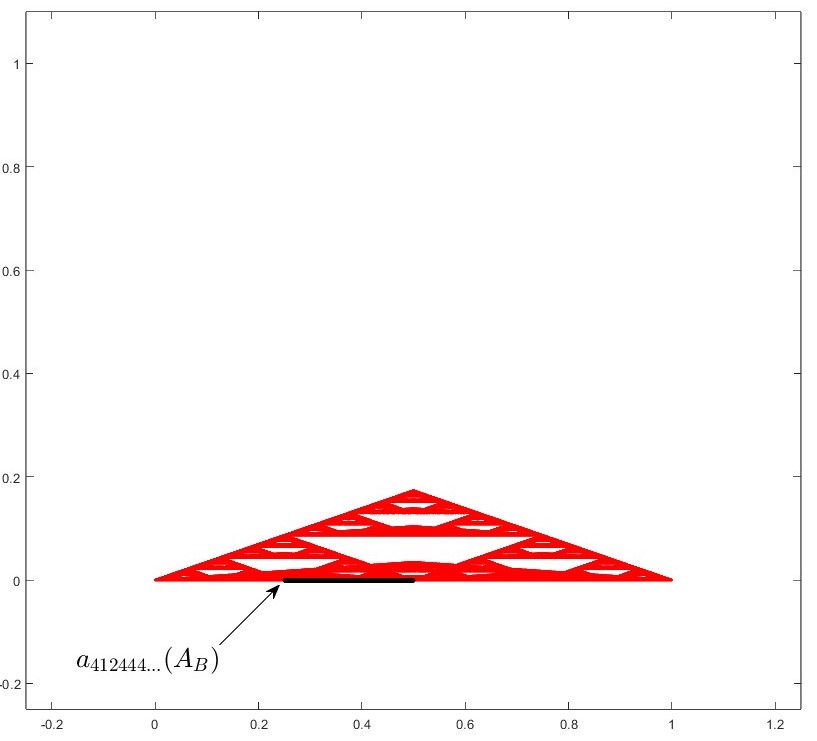}}
\subfloat[]{
\includegraphics[width=0.33\textwidth]{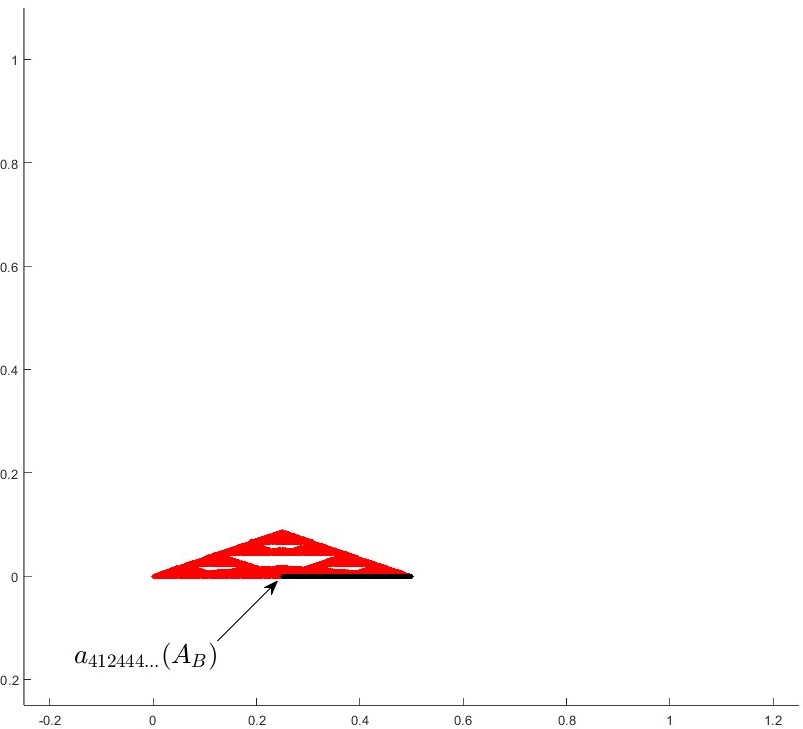}}
\subfloat[]{
\includegraphics[width=0.33\textwidth]{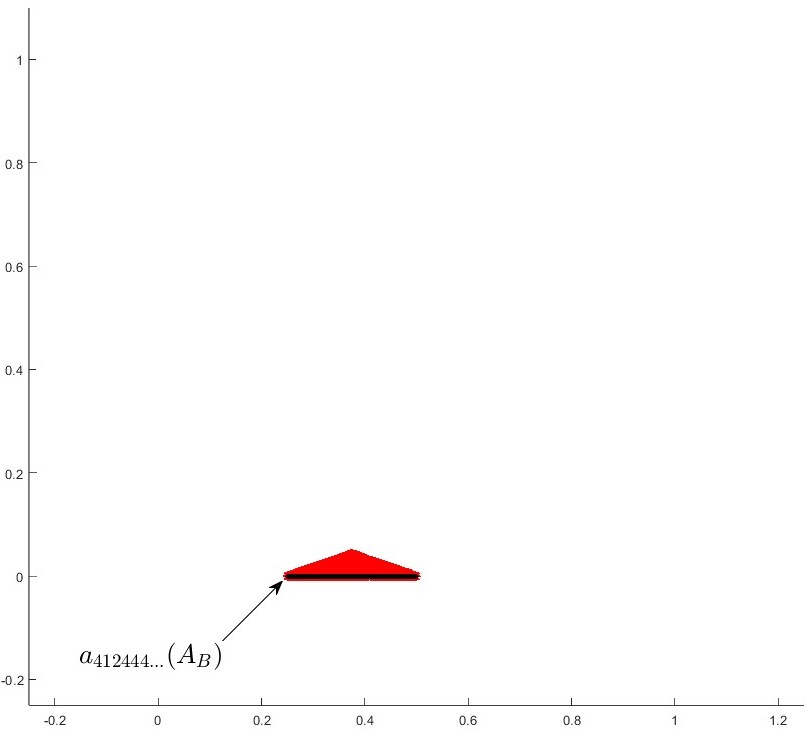}}\\
\subfloat[]{
\includegraphics[width=0.33\textwidth]{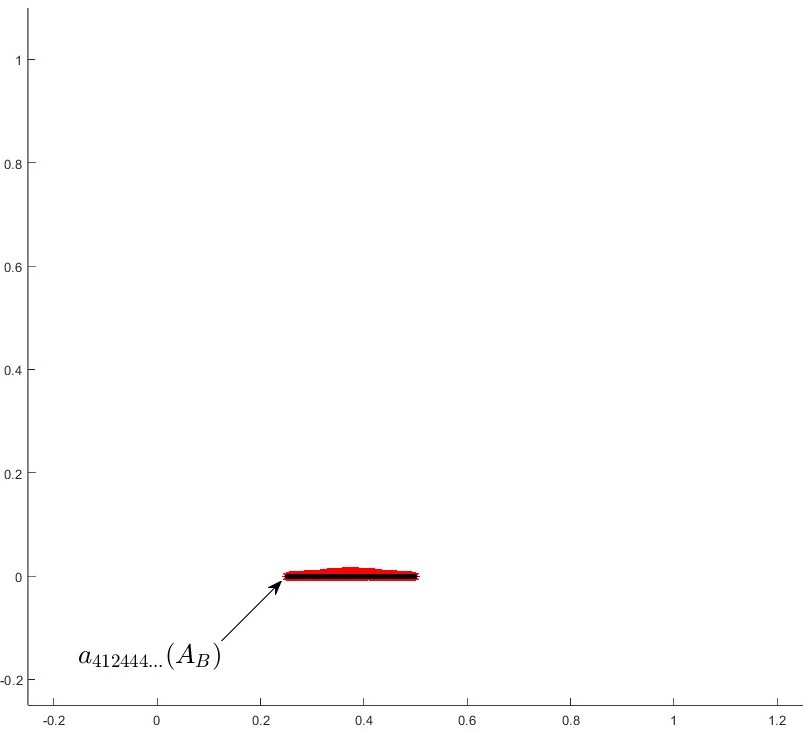}}
\subfloat[]{
\includegraphics[width=0.33\textwidth]{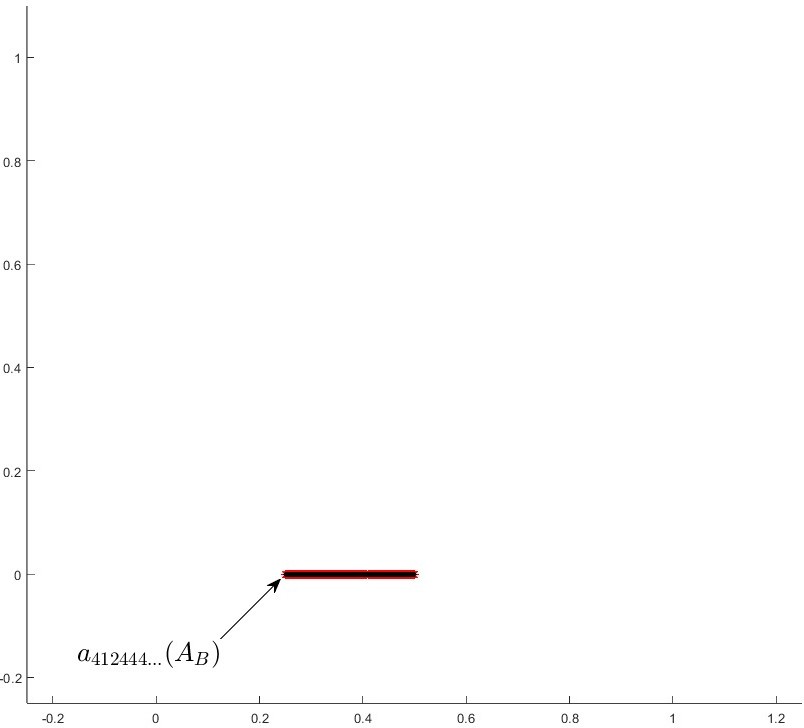}}
\caption{}
\end{figure}
\textbf{D}. Finally we note that the concatenation of the images of the
approximates of $a_{\alpha }((0,0))$ for 50000 randomly generated elements $%
\alpha \in \Sigma _{1}(I,J)$ with the images of $\mathcal{A}_{\alpha
}((0,0)) $ for another 50000 randomly generated elements $\alpha \in \Sigma
_{0}(I,J)$ (see Figure 5(a)) looks similar to the graphical representation of $%
F_{\mathcal{S}}^{[20]}(\{(0,0)\})$ (see Figure 5(b)). This remark is in
accordance with Theorem 3.12.
\vspace{-1em}
\begin{figure}[H]
    \centering
\subfloat[]{
\includegraphics[width=0.30\textwidth]{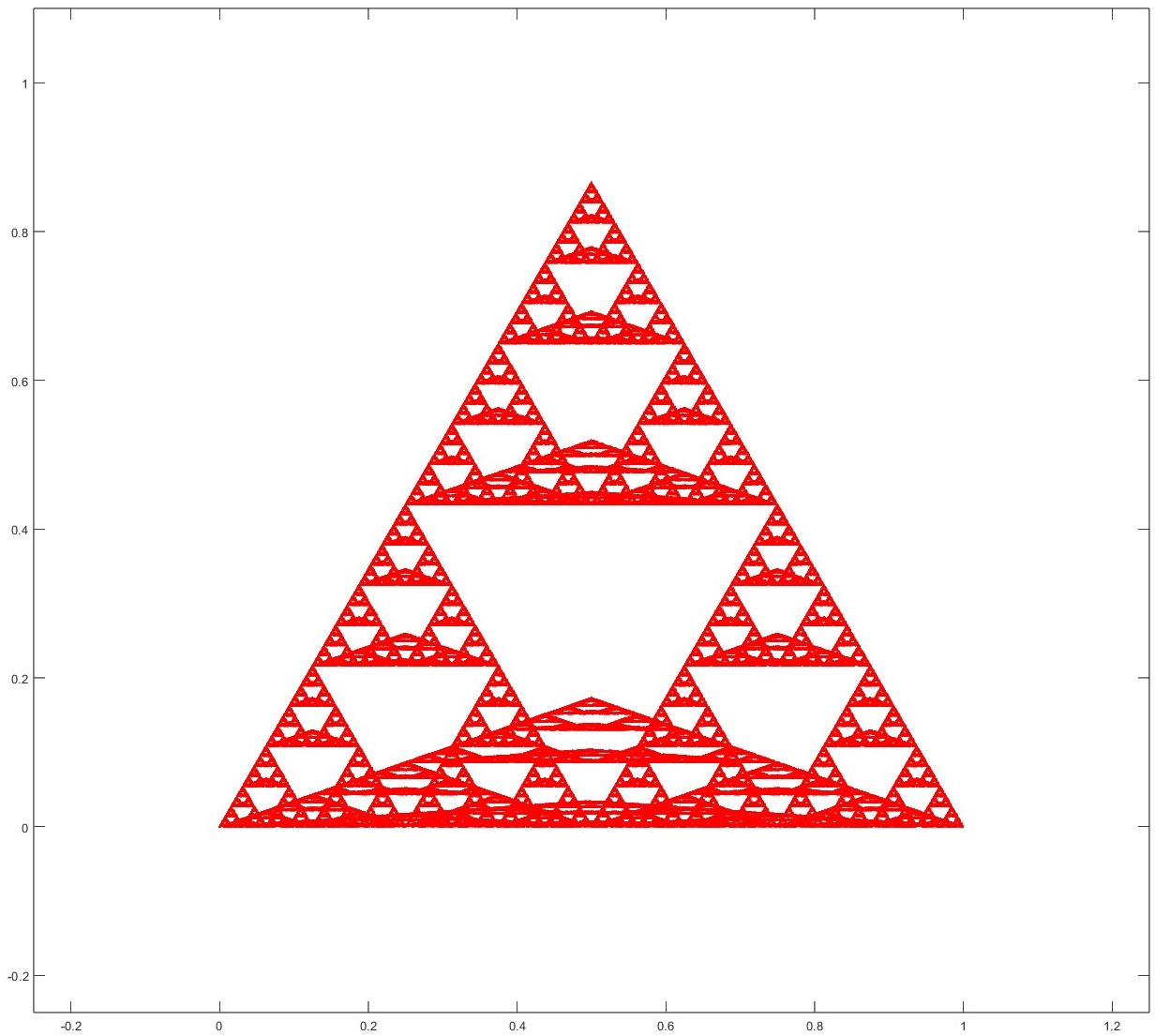}}
\subfloat[]{
\includegraphics[width=0.30\textwidth]{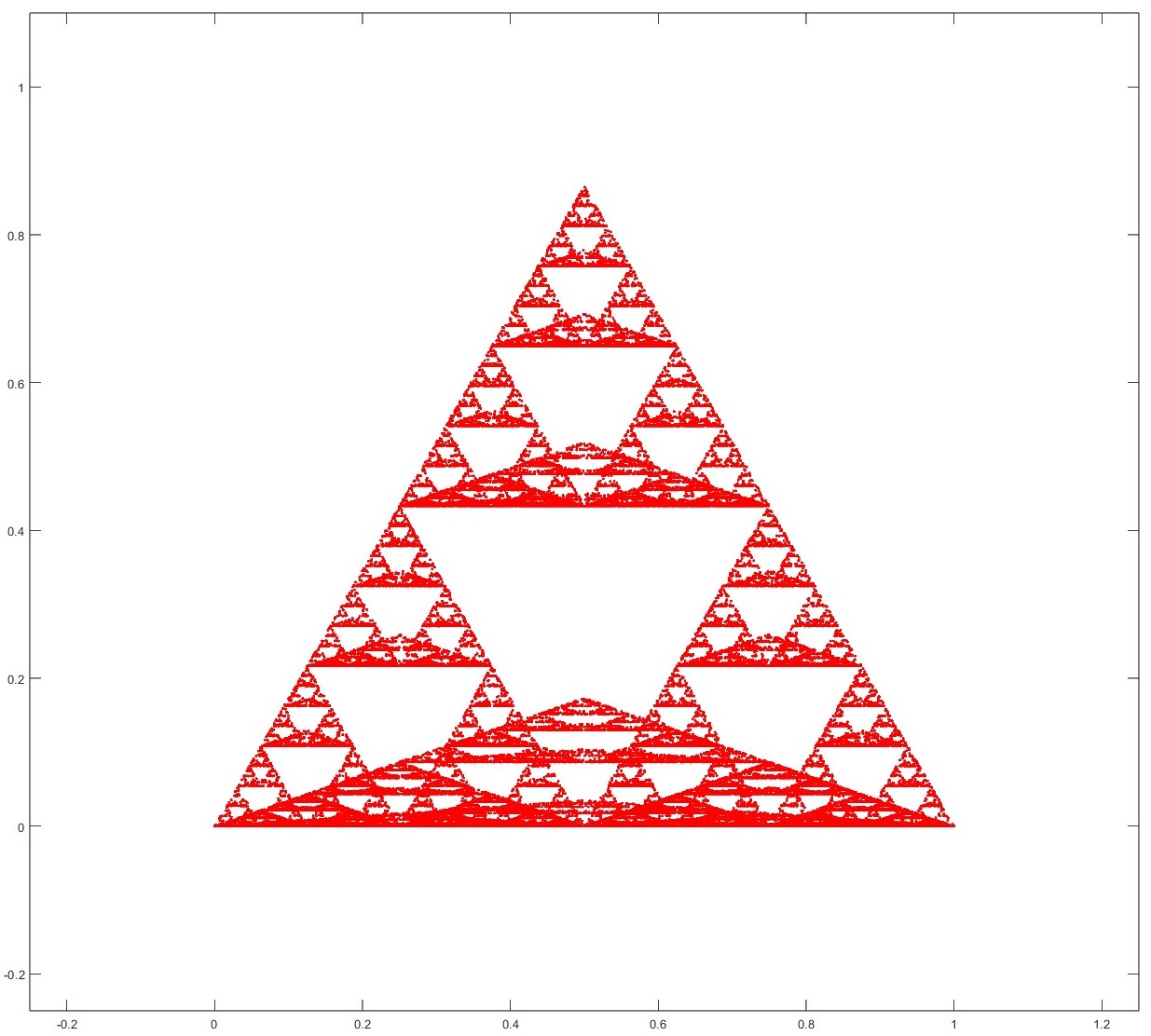}}
\caption{}
\end{figure}



\bigskip
\textbf{Statements and Declarations}

\textbf{Funding}: The authors declare that no funds, grants, or other support were received during the preparation of this manuscript.

\textbf{Competing Interests}: The authors declare no competing interests.

\textbf{Data availability}: Not applicable.

\textbf{References}

\bigskip

[1] B. Anghelina, R. Miculescu, A. Mihail, On the fractal operator of a
mixed possibly infinite iterated function system, Rev. Real Acad. Cienc.
Exactas Fis. Nat. Ser. A-Mat., \textbf{119} (2025), 31.

[2] H. Fernau, Infinite iterated function systems, Math. Nachr., \textbf{170}
(1994), 79--91.

[3] G. Gw\'{o}\'{z}d\'{z}- \L ukowska, J. Jachymski, The Hutchinson-Barnsley
theory for infinite iterated function systems, Bull. Aust. Math. Soc., 
\textbf{72} (2005), 441--454.

[4] M. Hille, Remarks on limits sets of infinite iterated functions systems,
Monatsh. Math., \textbf{168} (2012), 215--237.

[5] J. Hutchinson, Fractals and self similarity, Indiana Univ. Math. J., 
\textbf{30} (1981), 713--747.

[6] M. Iosifescu, Iterated function systems. A critical survey, Math. Rep.
(Bucur.), \textbf{11} (2009), 181--229.

[7] A. Kameyama, Distances on topological self-similar sets and the kneading
determinants, J. Math. Kyoto Univ., \textbf{40} (2000), 601-672.

[8] B. Kieninger, Iterated Function Systems on Compact Hausdorff Spaces, Ph.
D. thesis, University of Augsburg, Aachen: Shaker-Verlag, 2002.

[9] H. Kunze, D. La Torre, F. Mendivil, E. Vrscay, Fractal-Based Methods in
Analysis, Springer, 2012.

[10] K. Le\'{s}niak, Infinite iterated function systems: a multivalued
approach, Bull. Pol. Acad. Sci., Math., \textbf{52} (2004), 1--8.

[11] K. Le\'{s}niak, N. Snigireva, F. Strobin, Weakly contractive iterated
function systems and beyond: a manual, J. Differ. Equ. Appl., \textbf{26}
(2020), 1114--1173.

[12] D. Mauldin, M. Urba\'{n}ski, Dimensions and measures in infinite
iterated function systems, Proc. Lond. Math. Soc., \textbf{73} (1996),
105--154.

[13] D. Mauldin, M. Urba\'{n}ski, Graph directed Markov Systems, Cambridge
University Press, 2003.

[14] F. Mendivil, A generalization of IFS with probabilities to infinitely
many maps, Rocky Mountain J. Math., \textbf{28} (1998), 1043--1051.

[15] A. Mihail, R. Miculescu, The shift space for an infinite iterated
function system, Math. Rep. Bucur., \textbf{61} (2009), 21-32.

[16] A. Mihail, I. Savu, $\varphi $-contractive parent-child possibly
infinite IFSs and orbital $\varphi $-contractive possibly infinite IFSs,
Fixed Point Theory, \textbf{25} (2024), 229-248.

[17] F. Strobin, Contractive iterated function systems enriched with
nonexpansive maps, Result. Math., \textbf{76} (2021), 153.

\bigskip

Bogdan Cristian Anghelina

Faculty of Mathematics and Computer Science

Transilvania University of Bra\c{s}ov

Iuliu Maniu Street, nr. 50, 500091, Bra\c{s}ov, Romania

E-mail: bogdan.anghelina@unitbv.ro

\bigskip

Radu Miculescu

Faculty of Mathematics and Computer Science

Transilvania University of Bra\c{s}ov

Iuliu Maniu Street, nr. 50, 500091, Bra\c{s}ov, Romania

E-mail: radu.miculescu@unitbv.ro

\bigskip

Alexandru Mihail

Bucharest\ University

Faculty of Mathematics and Computer Science

Str. Academiei\ 14

010014 Bucharest, Romania

E-mail: mihail\_alex@yahoo.com

\end{document}